\documentclass{article}

\RequirePackage[OT1]{fontenc}
\RequirePackage{amsthm,amsmath,bm}
\RequirePackage{natbib}
\RequirePackage[colorlinks,citecolor=blue,urlcolor=blue]{hyperref}

\usepackage{amsmath,amssymb, a4wide}
\usepackage{amssymb,latexsym}
\usepackage{dsfont, color}
\usepackage{verbatim}
\setcitestyle{authoryear,open={(},close={)}} 
\usepackage{cancel}
\usepackage{appendix}
\usepackage{graphicx}
\usepackage{authblk}

\numberwithin{equation}{section}
\theoremstyle{plain}
\newtheorem{theorem}{Theorem}
\newtheorem{corollary}{Corollary}
\newtheorem{lemma}{Lemma}
\newtheorem{definition}{Definition}
\newtheorem{proposition}{Proposition}
\newtheorem{remark}{Remark}[section]

\newcommand{\R}{\mathbb{R}}

\newcommand{\dd}{\mathrm{d}}
\newcommand{\vc}{\mathrm{vec}}
\newcommand{\vch}{\mathrm{vech}}

\newcommand{\bbeta}{\bm\beta}
\newcommand{\beps}{\bm\epsilon}
\newcommand{\btheta}{\bm\theta}
\newcommand{\bxi}{\bm\xi}

\newcommand{\bgamma}{\bm\gamma}
\newcommand{\bmu}{\bm\mu}
\newcommand{\bSigma}{\mathbf{\Sigma}}

\newcommand{\by}{\mathbf{y}}
\newcommand{\bx}{\mathbf{x}}
\newcommand{\bs}{\mathbf{s}}
\newcommand{\bu}{\mathbf{u}}

\newcommand{\bt}{\mathbf{t}}
\newcommand{\bz}{\mathbf{z}}

\newcommand{\bb}{\mathbf{b}}

\newcommand{\bX}{\mathbf{X}}
\newcommand{\bV}{\mathbf{V}}

\newcommand{\bZ}{\mathbf{Z}}
\newcommand{\bI}{\mathbf{I}}
\newcommand{\bA}{\mathbf{A}}

\newcommand{\bC}{\mathbf{C}}
\newcommand{\bB}{\mathbf{B}}
\newcommand{\bD}{\mathbf{D}}

\newcommand{\bT}{\mathbf{T}}

\newcommand{\E}{\mathbb{E}}

\title{Highly Efficient Estimators with High Breakdown Point for Linear Models with Structured Covariance Matrices}
\author[1]{Hendrik Paul Lopuha\"a}
\affil[1]{\emph{Delft University of Technology}}
\date{\today}

\begin{document}
\maketitle
\begin{abstract}
A unified approach is provided for a method of estimation of the regression parameter
in balanced linear models with a structured covariance matrix
that combines a high breakdown point with high asymptotic efficiency at models with multivariate normal errors.
Of main interest are linear mixed effects models, but our approach also includes
several other standard multivariate models, such as multiple regression,
multivariate regression, and multivariate location and scatter.
Sufficient conditions are provided for the existence of the estimators and corresponding functionals,
strong consistency and asymptotic normality is established,
and robustness properties are derived in terms of breakdown point and influence function.
All the results are obtained for general identifiable covariance structures and are established under mild conditions
on the distribution of the observations, which goes far beyond models with elliptically contoured densities.
Some results are new and others are more general than existing ones in the literature.
In this way, results on high breakdown estimation with high efficiency in a wide variety of multivariate models are completed and improved.
\end{abstract}

\section{Introduction}
Linear models are widely used and provide a versatile approach for analyzing correlated responses,
such as longitudinal data, growth data or repeated measurements.
In such models, each subject~$i$, $i=1,\ldots,n$, is observed at $k_i$ occasions,
and the vector of responses~$\by_i$ is assumed to arise from the model
\[
\by_i=\bX_i\bbeta+\bu_i,
\]
where $\bX_i$ is the design matrix for the $i$th subject and $\bu_i$ is a vector
whose covariance matrix can be used to model the correlation between the responses.
One possibility is the linear mixed effects model, in which the random effects together with the measurement error
yields a specific covariance structure depending on a vector $\btheta$ consisting of some unknown covariance parameters.
Other covariance structures may arise, for example if the $\bu_i$ are the outcome of a time series.
See e.g.,~\citet{jennrich&schluchter1986} or~\cite{fitzmaurice-laird-ware2011}, for several possible covariance structures.

Maximum likelihood estimation of $\bbeta$ and $\btheta$ has been studied, e.g., in~\cite{hartley&rao1967,rao1972,laird&ware1982},
see also~\cite{fitzmaurice-laird-ware2011,demidenko2013}.
To be resistant against outliers, robust methods have been investigated for the linear mixed effects models, e.g.,
in~\cite{pinheiro-liu-wu2001,copt2006high,copt&heritier2007,heritier-cantoni-copt-victoriafeser2009,koller2013robust,chervoneva2014,agostinelli2016composite}.
This mostly concerns S-estimators, originally introduced in the multiple regression context by~\cite{rousseeuw-yohai1984}
and extended to multivariate location and scatter in~\cite{davies1987,lopuhaa1989,fishbone2021}, to multivariate linear regression in~\cite{vanaelst&willems2005},
and to linear mixed effects models
in~\cite{copt2006high,heritier-cantoni-copt-victoriafeser2009,chervoneva2011,chervoneva2014}.
A unified approach to S-estimation in balanced linear models with structured covariances can be found in~\cite{lopuhaa-gares-ruizgazenARXIVE2022}.

S-estimators are well known smooth versions of the minimum volume ellipsoid estimator~\cite{rousseeuw1985} that
are highly resistant against outliers and are asymptotically normal at $\sqrt{n}$-rate.
Unfortunately, the choice of the tuning constant corresponding to an S-estimator, forces a trade-off between robustness and efficiency.
For this reason,
remedies have been developed that retain the high breakdown point of the S-estimator
and improve the efficiency of the regression estimator in a second step.
One possibility are MM-estimators, introduced by~\cite{yohai1987} in the multiple regression setup.
Extensions to multivariate location and scatter can be found
in~\cite{lopuhaa1992highly,tatsuoka&tyler2000,SalibianBarrera-VanAelst-Willems2006,fishbone2021}.
An extension to linear mixed effects models was discussed in~\cite{copt&heritier2007}
and to multivariate linear regression by~\cite{kudraszow-maronna2011}.
An application of MM-estimation to emitter localization can be found in~\cite{parkchang2021}.

We will extend the approaches in~\cite{lopuhaa1992highly} and~\cite{copt&heritier2007} to balanced linear models with structured covariance matrices,
and postpone MM-estimation for unbalanced models to a future manuscript.
The balanced setup is already quite flexible and includes several specific multivariate statistical models.
Of main interest are high breakdown estimators with high normal efficiency for linear mixed effects models,
but our approach also includes high breakdown estimators in several other standard multivariate models, such as multiple regression,
multivariate linear regression, and multivariate location and scatter.
We provide sufficient conditions for the existence of the estimators and corresponding functionals,
establish their asymptotic properties, such as consistency and asymptotic normality,
and derive their robustness properties in terms of breakdown point and influence function.
All results are obtained for a large class of identifiable covariance structures, and are established under very mild conditions
on the distribution of the observations, which goes far beyond models with elliptically contoured densities.
In this way, some of our results are new and others are more general than existing ones in the literature.

The paper is organized as follows.
In Section~\ref{sec:structured covariance model}, we explain the model in detail and
provide some examples of standard multivariate models that are included in our setup.
In Section~\ref{subsec:def MM Yohai} we define the regression M-estimator and M-functional
and in Section~\ref{subsec:Existence MM Yohai} we give conditions under which they exist.
In Section~\ref{subsec:continuity general} we establish continuity of the regression M-functional,
which is then used to obtain consistency of the regression M-estimator.
Section~\ref{subsec:BDP Yohai} deals with the breakdown point.
Section~\ref{sec:score equations} provides the preparation for
Sections~\ref{subsec:IF MM Yohai} and~\ref{sec:asymptotic normality},
in which we determine the influence function and establish asymptotic normality.
Finally, in Section~\ref{sec:application},
we investigate the performance of the estimators by means of a simulation
and an application to data from a trial on the treatment of lead-exposed children.

\section{Balanced linear models with structured covariances}
\label{sec:structured covariance model}
We consider independent observations $(\by_1,\bX_1),\ldots,(\by_n,\bX_n)$,
for which we assume the following model
\begin{equation}
\label{def:model}
\by_i
=
\bX_i\bbeta+\bu_i,
\quad
i=1,\ldots,n,
\end{equation}
where $\by_i\in\R^{k}$ contains repeated measurements for the $i$-th subject,
$\bbeta\in\R^q$ is an unknown parameter vector,
$\bX_i\in\R^{k\times q}$ is a known design matrix, and
the $\mathbf{u}_i\in\R^{k}$ are unobservable independent mean zero random vectors with
covariance matrix $\bV\in\text{PDS}(k)$,
the class of positive definite symmetric $k\times k$ matrices.
The model is balanced in the sense that all~$\mathbf{y}_i$ have the same dimension.
Furthermore, we consider a structured covariance matrix, that is,
the matrix $\bV=\bV(\btheta)$ is a known function of unknown covariance parameters combined in a vector $\btheta\in\R^l$.
We first discuss some examples that are covered by this setup.

An important case of interest is the (balanced) linear mixed effects model.
For a general formulation covered by our setup, see~\cite{lopuhaa-gares-ruizgazenARXIVE2022}.
A specific example is the model
\begin{equation}
\label{def:linear mixed effects model Copt}
\by_i=\bX_i\bbeta+\sum_{j=1}^r \bZ_j\gamma_{ij}+\beps_i,
\quad
i=1,\ldots,n,
\end{equation}
considered in~\cite{copt&heritier2007}.
This model arises from $\bu_i=\sum_{j=1}^r \bZ_j\gamma_{ij}+\beps_i$, for $i=1,\ldots,n$,
where the~$\bZ_j$'s are known $k\times g_j$ design matrices and
the $\gamma_{ij}\in\R^{g_j}$ are independent mean zero random variables with covariance matrix $\sigma_j^2\bI_{g_j}$,
for $j=1,\ldots,r$, independent from $\beps_i$, which has mean zero and covariance matrix~$\sigma_0^2\bI_k$.
In this case, $\bV(\btheta)=\sum_{j=1}^r\sigma_j^2\bZ_j\bZ_j^T+\sigma_0^2\bI_k$ and~$\btheta=(\sigma_0^2,\sigma_1^2,\ldots,\sigma_r^2)$.

Another example of~\eqref{def:model}, is the multivariate linear regression model
\begin{equation}
\label{def:multivariate linear regression model}
\by_i=\bB^T\bx_i+\bu_i,
\qquad
i=1,\ldots,n,
\end{equation}
considered in~\cite{kudraszow-maronna2011}, where $\bB\in\R^{q\times k}$ is a matrix of unknown parameters, $\bx_i\in\R^q$ is known,
and~$\mathbf{u}_i$, for $i=1,\ldots,n$, are independent mean zero random variables with
covariance matrix~$\bV(\btheta)=\bC\in\text{PDS}(k)$.
In this case, the vector of unknown covariance parameters is given  by
\begin{equation}
\label{def:theta for unstructured}
\btheta=\vch(\bC)=(c_{11},\ldots,c_{1k},c_{22},\ldots,c_{kk})^T\in\R^{\frac12k(k+1)},
\end{equation}
The model can be obtained as a special case of~\eqref{def:model}, by taking
$\bX_i=\bx_i^T\otimes \bI_k$ and $\bbeta=\vc(\bB^T)$, where $\vc(\cdot)$ is the $k^2$-vector that stacks the columns of a matrix.
Clearly, the multiple linear regression model considerd in~\cite{yohai1987} is a special case with $k=1$.

Also the multivariate location-scale model, as considered in~\cite{lopuhaa1992highly}
(see also~\cite{SalibianBarrera-VanAelst-Willems2006,tatsuoka&tyler2000}),
can be obtained as a special case of~\eqref{def:model},
by taking $\bX_i=\bI_k$, the $k\times k$ identity matrix.
In this case, $\bbeta\in\R^k$ is the unknown location parameter and covariance matrix
$\bV(\btheta)=\bC\in\text{PDS}(k)$, with $\btheta$ as in~\eqref{def:theta for unstructured}.

Model~\eqref{def:model} also includes examples, for which
$\bu_1,\ldots,\bu_n$ are generated by a time series.
An example is the case where $\bu_i$ has a covariance matrix with elements
$v_{st}=\sigma^2\rho^{|s-t|}$,
for $s,t=1,\ldots,n$.
This arises when the $\bu_i$'s are generated by an autoregressive process of order one.
The vector of unknown covariance parameters is $\btheta=(\sigma^2,\rho)\in(0,\infty)\times[-1,1]$.
A general stationary process leads to
$v_{st}=\theta_{|s-t|+1}$,
for $s,t=1,\ldots,n$,
in which case $\btheta=(\theta_1,\ldots,\theta_k)^T\in\R^k$,
where $\theta_{|s-t|+1}$ represents the autocovariance over lag~$|s-t|$.

Throughout the manuscript we will assume that the parameter $\btheta$ is identifiable in the sense that,
$\bV(\btheta_1)=\bV(\btheta_2)$ implies $\btheta_1=\btheta_2$.
This is true for all examples mentioned above.

\section{Definitions}
\label{subsec:def MM Yohai}
The definition of MM-estimators involves the use of a single real-valued function
$\rho$, or the use of multiple real-valued functions $\rho_0$ and $\rho_1$.
Moreover, depending on the specific statistical model of interest,
the breakdown behavior of the corresponding MM-estimator may depend
on whether the $\rho$-functions are bounded or unbounded.
Since we intend to include both possibilities, we first discuss them both.

\subsection{Bounded and unbounded $\rho$-functions}
\label{subsec:bounded rho}
\cite{yohai1987} defines the regression MM-estimator in multiple stages.
By means of a function $\rho_0$, an M-estimator of scale is determined from residuals,
that are obtained from an initial high breakdown regression estimator.
Given the M-estimator of scale, a final regression M-estimator is determined by means of a function~$\rho_1$.
The conditions imposed on the two $\rho$-functions are similar to the following conditions.
\begin{quote}
\begin{itemize}
\item[(R-BND)]
$\rho$ is symmetric around zero with $\rho(0)=0$ and $\rho$ is continuous at zero.
There exists a finite constant $c>0$, such that $\rho$ is strictly increasing on $[0,c]$ and constant on~$[c,\infty)$;
put $a=\sup\rho=\rho(c)$.
\end{itemize}
\end{quote}
In addition, the two $\rho$-functions are related.
Suitable tuning of the bounded function $\rho_0$ ensures a high breakdown point of the scale M-estimator,
and by imposing the relationship between $\rho_0$ and~$\rho_1$,
the final regression M-estimator inherits the high breakdown point from the scale M-estimator.
Typical choices for bounded $\rho_0$ and $\rho_1$ that satisfy (R-BND), can be determined from Tukey's biweight,
defined as
\begin{equation}\label{def:biweight}
\rho_{\mathrm{B}}(s;c)
=
\begin{cases}
\displaystyle{\frac{s^2}2-\frac{s^4}{2c^2}+\frac{s^6}{6c^4}}, & |s|\leq c,\\
\\[-10pt]
\dfrac{c^2}{6} & |s|>c,
\end{cases}
\end{equation}
by taking $\rho_0(d)=\rho_{\mathrm{B}}(d;c_0)$ and $\rho_1(d)=\rho_{\mathrm{B}}(d;c_1)$,
where the cut-off constants are chosen such that $0<c_0<c_1<\infty$.
The cut-off constant $c_0$ can be tuned such that the MM-estimator inherits the breakdown point of the initial
regression estimator, whereas the constant $c_1$ can be tuned such that the MM-estimator has high efficiency at the
model with Gaussian errors.

In~\cite{lopuhaa1992highly} this idea has been extended to multivariate location and scatter by
determining a location M-estimator after first obtaining a high breakdown covariance estimator.
After rescaling the observations with the initial covariance estimator, a location M-estimator is obtained by
minimizing an object function that involves only a single $\rho$-function that satisfies the following condition.
\begin{quote}
\begin{itemize}
\item[(R-UNB)]
$\rho$ is symmetric, $\rho(0)=0$ and $\rho(s)\to\infty$, as $s\to\infty$.
The functions $\rho'$ and $u(s)=\rho'(s)/s$ are continuous,
$\rho'\geq 0$ on $[0,\infty)$ and there exists a $s_0$ such that $\rho'$ is nondecreasing on $(0,s_0)$ and
nonincreasing on $(s_0,\infty)$.
\end{itemize}
\end{quote}
In view of the results found by~\cite{huber1984},
an unbounded $\rho$-function is used in~\cite{lopuhaa1992highly}
to avoid that the breakdown point of the location M-estimator depends on the configuration of the sample,
which is the case for bounded $\rho$-functions.
With an unbounded $\rho$-function, the location M-estimator is shown to inherit the breakdown point of the initial covariance estimator.
A typical choice of an unbounded $\rho$-function that satisfies~(R-UNB) is
\begin{equation}\label{def:huber psi}
\rho_{\mathrm{H}}(s;c)
=
\begin{cases}
\dfrac{s^2}2, & |s|\leq c,\\
\\[-10pt]
-\dfrac{c^2}2+c|s|, & |s|>c,
\end{cases}
\end{equation}
whose derivative $\psi_\mathrm{H}=\rho_\mathrm{H}'$ is a bounded monotone function known as Huber's $\psi$-function.
The constant $c$ can be tuned such that the location M-estimator has high efficiency at the multivariate normal distribution.

\cite{tatsuoka&tyler2000} and~\cite{SalibianBarrera-VanAelst-Willems2006},
propose a different version of location MM-estimators also using bounded $\rho$-functions.
Instead of using the entire covariance matrix as auxiliary statistic,
they estimate the shape of the scatter matrix along with location parameter
and only use a univariate auxiliary estimator for the scale of the scatter matrix.
In~\cite{SalibianBarrera-VanAelst-Willems2006} it is shown that the location and shape estimators
in the second step inherit the breakdown point of the
initial estimators used in the first step.
\cite{kudraszow-maronna2011} use a similar version for multivariate linear regression and also
establish that the regression and shape estimators in the second step inherit the breakdown point of the
initial estimators used in the first step.
\cite{copt&heritier2007} treat regression MM-estimators in the context of linear mixed effects models.
They allow both bounded and unbounded $\rho$-functions and briefly discuss the pros and cons,
but do not explicitly derive the breakdown point.

Extending the approach in~\cite{lopuhaa1992highly}  and~\cite{copt&heritier2007} to the regression parameter $\bbeta$ in the current setup~\eqref{def:model}
seems straightforward.
First obtain a high breakdown structured covariance estimator and determine a regression M-estimator from the re-scaled observations.
However, in order to make sure that the resulting M-estimator inherits the breakdown point from the initial covariance estimator,
the use of an unbounded $\rho$-function, as in~\cite{lopuhaa1992highly}, does not seem to be suitable.
The presence of the design matrices $\bX_i$ in the object function to be minimized, makes things more complex than
for multivariate location.
Alternatively, one could minimize a single object function based on a bounded $\rho$-function.
However, in view of the results in~\cite{huber1984},
in this case it seems difficult to ensure that the resulting regression M-estimator inherits the breakdown point from
the initial covariance estimator.

\subsection{The regression M-estimator and corresponding M-functional}
We start by representing our observations as points in $\R^k\times\R^{kq}$  in the following way.
For $r=1,\ldots, k$, let $\bx_r^T$ denote the $r$-th row of the $k\times q$ matrix $\bX$,
so that $\bx_r\in\R^q$.
We represent the pair $\mathbf{s}=(\mathbf{y},\bX)$ as an element in $\R^k\times\R^{kq}$ defined by
$\bs^T=(\by^T,  \bx_{1}^T,\ldots,  \bx_{k}^T)$.
In this way our observations can be represented as $\bs_1,\ldots,\bs_n$, with $\bs_i=(\by_i,\bX_i)\in\R^k\times\R^{kq}$.

We will show that an approach similar to~\cite{yohai1987}, using two bounded $\rho$-functions that are suitably related
turns out to be helpful.
For the moment, we intend to include both bounded as well as unbounded $\rho$-functions in our approach.
In order to do so, the estimator for $\bbeta$ is defined in two stages as follows.
\begin{definition}
\label{def:MM-estimator general}
Let $\bV_{0,n}$ be a (high breakdown) positive definite symmetric covariance estimator.
For a function $\rho_1:\R\to[0,\infty)$, define $\bbeta_{1,n}$ as the vector that minimizes
\begin{equation}
\label{eq:MM estimator general}
R_n(\bbeta)
=
\frac{1}{n}
\sum_{i=1}^{n}
\rho_1\left(
\sqrt{(\by_i-\bX_i\bbeta)^T\bV_{0,n}^{-1}(\by_i-\bX_i\bbeta)}
\right).
\end{equation}
\end{definition}
At this point, $\rho_1$ can be either bounded or unbounded.
Later on, we will further specify under what conditions on $\rho_1$, several properties hold for~$\bbeta_{1,n}$.
Note that one may choose any initial (high breakdown) covariance estimator, but in our setup we typically think of a structured
covariance estimator~$\bV_{0,n}=\bV(\btheta_{0,n})$,
where $\btheta_{0,n}$ is an initial estimator for the vector of covariance parameters.
This means that $\bV_{0,n}$ is not necessarily affine equivariant, and similarly for~$\bbeta_{1,n}$.
However, it is not difficult to see that $\bbeta_{1,n}$ is regression equivariant, i.e.,
\[
\bbeta_{1,n}(\{(\by_i+\bX_i\bb,\bX_i),i=1,\ldots,n\})=\bbeta_{1,n}(\{(\by_i,\bX_i),i=1,\ldots,n\})+\bb,
\]
for all $\bb\in\R^q$.
The corresponding functional is defined similarly.
\begin{definition}
\label{def:MM-functional general}
Let $\bV_0(P)$ be a positive definite symmetric covariance functional.
For a function $\rho_1:\R\to[0,\infty)$, define $\bbeta_1(P)$ as the vector that minimizes
\begin{equation}
\label{eq:MM-functional general}
R_{P}(\bbeta)
=
\int
\rho_1\left(
\sqrt{(\by-\bX\bbeta)^T\bV_0(P)^{-1}(\by-\bX\bbeta)}
\right)
\,\dd
P(\by,\bX).
\end{equation}
\end{definition}
The functional $\bbeta_1(P)$ is regression equivariant in the sense that
\[
\bbeta_1(P_{\by+\bX\bb,\bX})=\bbeta_1(P_{\by,\bX})+\bb,
\]
for all $\bb\in\R^q$,
where $P_{\by,\bX}$ denotes the distribution of $(\by,\bX)$.
Clearly, if one takes $P=\mathbb{P}_n$, the empirical measure of
the sample $(\by_1,\bX_1),\ldots,(\by_n,\bX_n)$, then $\bbeta_1(\mathbb{P}_n)=\bbeta_{1,n}$.

As before, $\bV_0(P)$ can be any covariance functional,
but in our setup we typically think of a structured covariance
$\bV_0(P)=\bV(\btheta_0(P))$, where
$\btheta_0(P)$ is an initial functional representing
the vector of covariance parameters.
An example of an estimator $\btheta_{0,n}$ and corresponding functional~$\btheta_0(P)$
that yield a high breakdown structured covariance estimator~$\bV(\btheta_{0,n})$,
is the S-estimator $\btheta_{0,n}$ and its corresponding functional proposed in~\cite{lopuhaa-gares-ruizgazenARXIVE2022}.

Definitions~\ref{def:MM-estimator general} and~\ref{def:MM-functional general} coincide with the ones for the multivariate location M-estimator in~\cite{lopuhaa1992highly},
when we choose $\bX_i=\bI_k$ and $\bV(\btheta)=\bC\in\text{PDS}(k)$, with $\btheta$ as in~\eqref{def:theta for unstructured}.
For the multiple linear regression model~\eqref{def:multivariate linear regression model},
it follows that if $\bbeta_{1,n}$ exists, then it satisfies score equation~(2.6) and equation~(2.7) in~\cite{yohai1987}.
Similarly, for the linear mixed effects model~\eqref{def:linear mixed effects model Copt},
it follows that if $\bbeta_{1,n}$ exists, then it satisfies a score equation similar to equation~(8) in~\cite{copt&heritier2007}.

We should emphasize that score equations like (2.6) in~\cite{yohai1987} and (8) in~\cite{copt&heritier2007} are useful to
obtain asymptotic properties, but they do not guarantee that $\bbeta_{1,n}$ inherits the breakdown point of the estimators
used in the first step.
Breakdown behavior is typically established from the minimization problem in Definition~\ref{def:MM-estimator general} itself.
If this minimization problem has a solution $\bbeta_{1,n}$ and if this solution inherits the high breakdown point from $\bV_{0,n}$,
then~$\bbeta_{1,n}$ will be a zero of the corresponding score equation with a high breakdown point.
But just being a zero of the score equation does not ensure a high breakdown point.
Indeed, the breakdown point of the MM-estimators in the multiple linear regression model~\cite{yohai1987}
and the multivariate location-scale model~\cite{lopuhaa1992highly},
have been obtained from the respective minimization problems, and similarly
for the MM-estimators in~\cite{SalibianBarrera-VanAelst-Willems2006}
and~\cite{kudraszow-maronna2011}.
For the MM-estimator in the linear mixed effects model~\cite{copt&heritier2007}, the robustness properties have not been investigated.
In view of the fact that the use of bounded or unbounded $\rho$-functions may lead to different breakdown behavior,
the breakdown point of MM-estimators for the linear mixed effects model considered in~\cite{copt&heritier2007}
will be investigated in Section~\ref{subsec:BDP Yohai}.

\section{Existence}
\label{subsec:Existence MM Yohai}
Consider the functional $\beta_1(P)$, as defined in Definition~\ref{def:MM-functional general}.
We will establish existence of~$\beta_1(P)$,
where we allow both bounded and unbounded $\rho_1$.
Existence of the corresponding estimator $\bbeta_{1,n}$ will follow from this.
Also of interest is the special case in which $P$ is such that $\by\mid\bX$ has an elliptically contoured density of the form
\begin{equation}
\label{eq:elliptical}
f_{\bmu,\bSigma}(\by)
=
\text{det}(\bSigma)^{-1/2}
h\left(
(\by-\bmu)^T
\bSigma^{-1}
(\by-\bmu)
\right),
\end{equation}
with $\bmu=\bX\bbeta\in\R^k$ and $\bSigma=\bV(\btheta)\in\text{PDS}(k)$, and $h:[0,\infty)\to[0,\infty)$.
For the linear mixed effects model in~\cite{copt&heritier2007}, it is assumed that
$\by\mid\bX$ has a multivariate normal distribution, which is a special case of~\eqref{eq:elliptical}
with $h(t)=(2\pi)^{-k/2}\exp(-t/2)$.

For bounded $\rho_1$, we want to rule out the pathological case, where $P$ has all
of its mass outside the ellipsoid centered around the origin
with covariance structure $\bV_0(P)$ and radius~$c_1$.
To this end we require the following condition on~$P$.
\begin{quote}
\begin{itemize}
\item[(A)]
Suppose that
\[
R_P(\textbf{0})
=
\int
\rho_1\left(
\sqrt{\by^T\bV_0(P)^{-1}\by}
\right)\,\dd
P(\by,\bX)
<
\sup\rho_1.
\]
\end{itemize}
\end{quote}
Clearly, if $\by\mid\bX\sim \mathcal{N}(\bX\bbeta,\bSigma)$ and $\rho_1$ satisfies (R-BND),
this condition is trivially fulfilled.
We then have the following theorem for bounded~$\rho_1$.
\begin{theorem}
\label{th:existence MM bounded rho}
Let $\rho_1:\R\to[0,\infty)$ satisfy condition (R-BND),
and suppose that $\bX$ has full rank with probability one.
\begin{itemize}
\item[(i)]
If  $P$ satisfies~(A), then there is at least one  vector $\bbeta_1(P)$ that minimizes $R_P(\bbeta)$.
\item[(ii)]
When $P$ is such that $\by\mid\bX$ has an elliptically contoured density from~\eqref{eq:elliptical} with parameters $\bmu=\bX\bbeta$ and $\bSigma$,
and if $\bV_0(P)=\bSigma$, then $R_P(\bb)\geq R_P(\bbeta)$, for all~$\bb\in\R^q$.
When $h$ in~\eqref{eq:elliptical} and $-\rho_1$ have a common point of decrease, then $R_P(\bb)$ is uniquely minimized
by $\bbeta_1(P)=\bbeta$.
\end{itemize}
\end{theorem}
\begin{proof}
(i)
Let $0<\lambda_1<\infty$ be the largest eigenvalue of $\bV_0(P)$,
and let $\lambda_k(\bX^T\bX)>0$ denote the smallest eigenvalue of $\bX^T\bX$.
Let $\|\cdot\|$ denote the Euclidean norm.
Then we have that
\begin{equation}
\label{eq:lower bound d Yohai}
\begin{split}
\sqrt{(\by-\bX\bb)^T\bV_0(P)^{-1}(\by-\bX\bb)}
&\geq
\frac{\|\by-\bX\bb\|}{\sqrt{\lambda_1}}
\geq
\frac{\|\bX\bb\|-\|\by\|}{\sqrt{\lambda_1}}\\
&\geq
\frac{1}{\sqrt{\lambda_1}}
\left(
\|\bb\|\sqrt{\lambda_k(\bX^T\bX)}-\|\by\|
\right).
\end{split}
\end{equation}
Then by dominated convergence and (R-BND), it follows that
\begin{equation}
\label{eq:limit RP}
\lim_{\|\bb\|\to\infty}
R_P(\bb)
=
\int
\lim_{\|\bb\|\to\infty}
\rho_1\left(
\sqrt{(\by-\bX\bb)^T\bV_0(P)^{-1}(\by-\bX\bb)}
\right)
\,\dd
P(\by,\bX)\\
=
\sup\rho_1.
\end{equation}
According to condition~(A),
this means that there exists a constant $M>0$, such that
\begin{equation}
\label{eq:def M Yohai}
R_P(\bb)
>
R_P(\mathbf{0}),
\quad
\text{for all }\|\bb\|>M.
\end{equation}
Therefore, for minimizing~$R_P(\bb)$ we may restrict ourselves to the set
$K=\{\bb\in\R^q:\|\bb\|\leq M\}$.
By dominated convergence and (R-BND), it also follows that $R_P(\bb)$ is continuous
on the compact set~$K$,
and therefore it must attain at least one minimum~$\bbeta_1(P)$.

(ii)
Write
\[
R_P(\bb)
=
\E_\bX
\left[
\E_{\by\mid\bX}
\left[
\rho_1\left(
\sqrt{(\by-\bX\bb)^T\bSigma^{-1}(\by-\bX\bb)}
\right)
\right]
\right].
\]
By change of variables $\by=\bSigma^{1/2}\bz+\bmu$, the inner conditional expectation can be written as
\[
\int
\rho_1\left(
\|\bSigma^{-1/2}(\by-\bX\bb)\|
\right)
f_{\bmu,\bSigma}(\by)\,\dd\by
=
\int
\rho_1\left(
\|\bz-\bSigma^{-1/2}\bX(\bb-\bbeta)\|
\right)
h(\bz^T\bz)\,\dd\bz.
\]
Next, we apply Lemma~4 from~\cite{davies1987} to the functions $\xi(d)=1-\rho_1(\sqrt{d})/a_0$ and $g=h$
and taking $\Lambda=\bI_k$.
Since $-\rho_1$ and $h$ have a common point of decrease, for all $\bX$, it follows that
\[
\int
\rho_1\left(
\|\bz-\bSigma^{-1/2}\bX(\bb-\bbeta)\|
\right)
h(\bz^T\bz)\,\dd\bz
\leq
\int
\rho_1\left(
\|\bz\|
\right)
h(\bz^T\bz)\,\dd\bz,
\]
with a strict inequality unless $\bSigma^{-1/2}\bX(\bb-\bbeta)=\mathbf{0}$,
i.e., unless $\bb=\bbeta$, since $\bX$ has full rank with probability one.
Finally, with the same change of variables $\bz=\bSigma^{-1/2}(\by-\bmu)$,
the right hand side can be written as
\[
\int
\rho_1\left(
\|\bz\|
\right)
h(\bz^T\bz)\,\dd\bz
=
\int
\rho_1\left(
\|\bSigma^{-1/2}(\by-\bX\bbeta)\|
\right)
f_{\bmu,\bSigma}(\by)\,\dd\by.
\]
After taking expectations $\E_\bX$, we conclude that $R_P(\bb)\leq R_P(\bbeta)$,
with a strict inequality, unless~$\bb=\bbeta$.
This proves the theorem.
\end{proof}
For bounded $\rho_1$, the function $R_P(\bbeta)$ in~\eqref{eq:MM-functional general} is well defined.
This is not necessarily true for unbounded $\rho_1$.
However, this will be the case when $P$ has a first moment.
For unbounded $\rho_1$ we have the following result.
\begin{theorem}
\label{th:existence MM unbounded rho}
Let $\rho_1:\R\to[0,\infty)$ satisfy condition (R-UNB).
Suppose that $\E_P\|\bs\|<\infty$ and that~$\bX$ has full rank with probability one.
\begin{itemize}
\item[(i)]
For every $\bbeta\in\R^q$ fixed, $R_P(\bbeta)<\infty$.
\item[(ii)]
There is at least one vector $\bbeta_1(P)$ that minimizes $R_P(\bbeta)$.
When $\rho_1$ is also strictly convex, then $\bbeta_1(P)$ is uniquely defined.
\item[(iii)]
When $P$ is such that $\by\mid\bX$ has an elliptically contoured density from~\eqref{eq:elliptical} with parameters $\bmu=\bX\bbeta$ and $\bSigma$,
and if $\bV_0(P)=\bSigma$, then $R_P(\bb)\geq R_P(\bbeta)$, for all $\bb\in\R^q$.
When $h$ in~\eqref{eq:elliptical} is strictly decreasing, then $R_P(\bb)$ is uniquely minimized
by $\bbeta_1(P)=\bbeta$.
\end{itemize}
\end{theorem}

\begin{proof}
Let $0<\lambda_k\leq \lambda_1<\infty$ be the smallest and largest eigenvalue of $\bV_0(P)$.

(i)
Condition (R-UNB) implies that $\rho_1(s)\leq\rho_1(s_0)$, for $s\in[0,s_0]$, and
that for $s>s_0$,
\begin{equation}
\label{eq:prop rho general}
\rho_1(s)=\int_{0}^{s_0}\rho_1'(t)\,\dd t+\int_{s_0}^{s}\rho_1'(t)\,\dd t
\leq
\rho_1(s_0)+(s-s_0)\rho_1'(s_0).
\end{equation}
Hence, for $\|\bV_0(P)^{-1/2}(\by-\bX\bbeta)\|>s_0$,
we have that
\begin{equation}
\label{eq:bound rho1}
\begin{split}
\rho_1\left(\|\bV_0(P)^{-1/2}(\by-\bX\bbeta)\|\right)
&\leq
\rho_1(s_0)+\|\by-\bX\bbeta\|\lambda_k^{-1/2}\rho_1'(s_0)-s_0\rho_1'(s_0)\\
&\leq
\rho_1(s_0)+(\|\by\|+\|\bX\|\cdot\|\bbeta\|)\lambda_k^{-1/2}\rho_1'(s_0).
\end{split}
\end{equation}
Since $\E_P\|\bs\|<\infty$, we find that for any $\bbeta\in\R^q$ fixed,
\[
\int
\rho_1\left(
\|\bV_0(P)^{-1/2}(\by-\bX\bb)\|
\right)
\,\dd P(\bs)
\leq
\rho_1(s_0)+(\E_P\|\by\|+\|\bbeta\|\E_P\|\bX\|)\lambda_k^{-1/2}\rho_1'(s_0)<\infty,
\]
which proves part~(i).

(ii)
We first argue that for minimizing $R_P(\bbeta)$,
we can restrict ourselves to a compact set.
Note that $R_P(\textbf{0})<\infty$, according to part~(i).
Now, suppose that $\|\bbeta\|>M$.
Then from~\eqref{eq:lower bound d Yohai},
\begin{equation}\label{eq:lower bound quadratic form}
\sqrt{(\by-\bX\bb)^T\bV_0(P)^{-1}(\by-\bX\bb)}
\geq
\frac{1}{\sqrt{\lambda_1}}
\left(
\|\bb\|\sqrt{\lambda_k(\bX^T\bX)}-\|\by\|
\right)
\geq
\frac{\sqrt{M}}{2\sqrt{\lambda_1}},
\end{equation}
on the set
\begin{equation}\label{def:AM}
A_M=
\left\{
(\by,\bX)\in\R^{k+qk}:
\lambda_k(\bX^T\bX)\geq 1/M;\,\|\by\|\leq \sqrt{M}/2
\right\}.
\end{equation}
Since $\bX$ has full rank with probability one,
$P(A_M)\to1$ and $\rho_1(\sqrt{M}/(2\sqrt{\lambda_1}))\to\infty$,
as $M\to\infty$,
according to (R-UNB).
This implies that for $M$ sufficiently large,
\[
R_P(\bbeta)
\geq
\rho_1(\sqrt{M}/(2\sqrt{\lambda_1}))P(A_M)
>
R_P(\textbf{0}).
\]
Therefore, that there exists a constant $M>0$, such that for minimizing $R_P(\bbeta)$ we may restrict ourselves to the compact set
$K=\{\bbeta\in\R^q:\|\bbeta\|\leq M\}$.
Since $\E_P\|\bs\|<\infty$, from~\eqref{eq:bound rho1} and dominated convergence,
it follows that
$R_P(\bbeta)$ is continuous on~$K$ and therefore it must attain at least one minimum $\bbeta_1(P)$ on the compact set $K$.
It is easily seen that strict convexity of $\rho_1$
implies strict convexity of $R_P$,
which means that $\bbeta_1(P)$ is unique.

(iii) Because $\bbeta_1(\cdot)$ is regression equivariant, we may assume that $\bbeta=\mathbf{0}$.
Write
\[
R_P(\bb)
=
\E_\bX
\left[
\E_{\by|\bX}
\left[
\rho_1\left(
\|\bV_0(P)^{-1/2}(\by-\bX\bb)\|
\right)
\right]\right].
\]
Since $\bz=\bV_0(P)^{-1/2}\by=\bSigma^{-1/2}\by$ has an elliptically contoured density with
parameters~$(\mathbf{0},\bI_k)$, the inner conditional expectation can be written as
\[
\iint
\big\{0\leq s\leq \rho_1(\|\bz-\bV_0(P)^{-1/2}\bX\bb\|)\big\}
f(\|\bz\|)
\,\dd s\,\dd\bz.
\]
From here on, we can copy the proof of Theorem~2.1 in~\cite{lopuhaa1992highly} and
conclude that
\begin{equation}
\label{eq:ineq RP general}
\E_{\by|\bX}
\left[
\rho_1\left(
\|\bV_0(P)^{-1/2}(\by-\bX\bb)\|
\right)
\right]
\geq
0,
\quad
\bX-\text{a.s.}
\end{equation}
It follows that $R_P(\bb)\geq 0=R_P(\mathbf{0})$.
When $f$ is strictly decreasing, similar to the proof of Theorem~2.1 in~\cite{lopuhaa1992highly},
it follows that inequality~\eqref{eq:ineq RP general} is strict, which yields $R_P(\bb)>0=R_P(\mathbf{0})$.
\end{proof}

A direct consequence of Theorems~\ref{th:existence MM bounded rho} and~\ref{th:existence MM unbounded rho}
is the existence of~$\bbeta_{1,n}$.
\begin{corollary}
\label{cor:existence MM estimator}
Let $(\by_1,\bX_1),\ldots,(\by_n,\bX_n)$ be a sample, such that $\bX_i$ has full rank for each $i=1,\ldots,n$.
\begin{enumerate}
\item[(i)]
If $\rho_1:\R\to[0,\infty)$ satisfies conditions (R-BND) and $R_n(\textbf{0})<\sup\rho_1$,
then there exists at least one $\bbeta_{1,n}$ that minimizes~$R_n(\bbeta)$.
\item[(ii)]
If $\rho_1:\R\to[0,\infty)$ satisfies condition (R-UNB), then there exists at least one $\bbeta_{1,n}$ that minimizes~$R_n(\bbeta)$.
When $\rho_1$ is also strictly convex, then $\bbeta_{1,n}$ is uniquely defined.
\end{enumerate}
\end{corollary}
\begin{proof}
Take $P$ equal to the empirical measure $\mathbb{P}_n$ of the sample $(\by_1,\bX_1),\ldots,(\by_n,\bX_n)$.
If~$\bX_i$ has full rank, for each $i=1,\ldots,n$, then if $\rho_1:\R\to[0,\infty)$ satisfies either (R-BND) or condition (R-UNB),
the corollary follows from Theorems~\ref{th:existence MM bounded rho}(i) and~\ref{th:existence MM unbounded rho}(ii).
\end{proof}
The condition $R_n(\mathbf{0})<\sup\rho_1$ is not very restrictive.
It rules out the pathological case of all observations being
outside the ellipsoid centered around the origin with covariance structure
$\bV_{0,n}$ and radius~$c_1$.

Existence of regression MM-estimators was not considered in~\cite{copt&heritier2007} for linear mixed effects models
or in~\cite{yohai1987} for multiple linear regression.
Their existence now follows from Corollary~\ref{cor:existence MM estimator}.
Existence of location MM-estimators, as obtained in~\cite{lopuhaa1992highly}, now also follows from
Corollary~\ref{cor:existence MM estimator} as a special case of part~(ii).
Existence of a slightly different MM-estimator has been established in~\cite{tatsuoka&tyler2000} for multivariate location,
and in~\cite{kudraszow-maronna2011} for multivariate linear regression.

\section{Continuity and Consistency}
\label{subsec:continuity general}
Consider a sequence $P_t$, $t\geq0$, of probability measures on $\R^k\times\R^{kq}$ that converges weakly to~$P$, as $t\to\infty$.
By continuity of the functional $\bbeta_1(P)$ we mean that $\bbeta_1(P_t)\to\bbeta_1(P)$, as $t\to\infty$.
An example of such a sequence is the sequence of empirical measures $\mathbb{P}_n$, $n=1,2,\ldots$, that converges weakly to $P$, almost surely.
Continuity of the functional~$\bbeta_1(P)$ for this sequence would then mean that the estimator $\bbeta_{1,n}$ is consistent,
i.e., $\bbeta_{1,n}=\bbeta_1(\mathbb{P}_n)\to\bbeta_1(P)$,
almost surely.
Furthermore, continuity of the $\bbeta_1(P)$ also provides a first step in deriving the influence function, in the sense that
$\bbeta_1(P_{\epsilon,\bs_0})\to\bbeta_1(P)$, as $\epsilon\downarrow0$, where
\begin{equation}
\label{def:perturbed P}
P_{\epsilon,\bs_0}=(1-\epsilon)P+\epsilon\delta_{\bs_0},
\end{equation}
with $\delta_{\bs_0}$
representing the Dirac measure at $\bs_0=(\by_0,\bX_0)$.

When $\rho_1$ is bounded, we can obtain continuity of the functional $\bbeta_1(P)$ for general weakly convergent sequences
$P_t$, $t\geq0$.
When $\rho_1$ is unbounded, this becomes more complicated, but we can still establish continuity for the sequence of
empirical measures $\mathbb{P}_n$, $n=1,2,\ldots$, and for the sequence $P_{\epsilon,\bs_0}$, for $\epsilon\downarrow0$.
For bounded~$\rho_1$ we have the following theorem.
\begin{theorem}
\label{th:continuity MM bounded rho}
Let $P_t$, $t\geq0$ be a sequence of probability measures on $\R^k\times\R^{kq}$ that converges weakly to~$P$, as $t\to\infty$.
Suppose that $\rho_1:\R\to[0,\infty)$ satisfies (R-BND) and
suppose that~$P$ is such that (A) holds and that $\bX$ has full rank with probability one.
Suppose that for $t$ sufficiently large, $\bV_0(P_t)$ exists and that
\begin{equation}
\label{eq:conv Vt Yohai}
\lim_{t\to\infty}\bV_0(P_t)=\bV_0(P).
\end{equation}
Then for $t$ sufficiently large, there exists at least one $\bbeta_1(P_t)$ that
minimizes $R_{P_t}(\bbeta)$.
If~$\bbeta_1(P)$ is the unique minimizer of $R_P(\bbeta)$, then for any sequence~$\bbeta_1(P_t)$, $t\geq 0$,
it holds that
\[
\lim_{t\to\infty}\bbeta_1(P_t)=\bbeta_1(P).
\]
\end{theorem}
\begin{proof}
Similar to Lemma~B.1 in~\cite{lopuhaa-gares-ruizgazenARXIVE2022}, one can show that
\begin{equation}
\label{eq:prop Lemma 2 Yohai}
\lim_{t\to\infty}
\int
\rho_1\left(
d(\bs,\bbeta_{t},\bV_{t})
\right)
\,\dd P_t(\bs)
=
\int
\rho_1\left(
d(\bs,\bbeta_{L},\bV_{L})
\right)
\,\dd P(\bs),
\end{equation}
for any sequence $(\bbeta_{t},\bV_{t})\to(\bbeta_{L},\bV_{L})$,
where
\begin{equation}
\label{def:mahalanobis}
d^2(\bs,\bbeta,\bV)=(\by-\bX\bbeta)^T\bV^{-1}(\by-\bX\bbeta).
\end{equation}
In particular, this yields that, for every $\bbeta\in\R^q$ fixed, it holds that
$R_{P_t}(\bbeta)\to R_P(\bbeta)$,
as $t\to\infty$.
We first show that there exists $M>0$, such that for minimizing $R_{P_t}(\bbeta)$,
we can restrict ourselves to $\|\bbeta\|\leq M$ for $t$ sufficiently large.
Consider the set $A_M$ defined in~\eqref{def:AM}.
Due to condition~(A) and the fact that $\bX$ has full rank, with probability one, for any $\eta>0$, we can find an $M>0$, such that
$\rho_1(\sqrt{M}/(2\sqrt{\lambda_1}))P(A_M)>R_P(\textbf{0})+2\eta$.
On the other hand, for any $\eta>0$, we have
\begin{equation}
\label{eq:bound difference}
|R_{P_t}(\mathbf{0})-R_P(\mathbf{0})|\leq \eta,
\end{equation}
for $t$ sufficiently large.
If $\bbeta$ minimizes $R_{P_t}(\bbeta)$, for $t$ sufficiently large, we must have $\|\bbeta\|\leq M$,
since otherwise, according to~\eqref{eq:lower bound quadratic form},
\[
R_{P_t}(\bbeta)
\geq
\rho_1(\sqrt{M}/(2\sqrt{\lambda_1}))P(A_M)
>
R_P(\textbf{0})+2\eta
\geq
R_{P_t}(\textbf{0})+\eta
>
R_{P_t}(\textbf{0}).
\]
Hence, for minimizing $R_{P_t}(\bbeta)$, we can restrict  to
the compact set $K=\{\bbeta\in\R^q:\|\bbeta\|\leq M\}$.
Furthermore, as in the proof of Theorem~\ref{th:existence MM bounded rho}, the function $R_{P_t}(\bbeta)$ is
continuous on the compact set~$K$, and must therefore attain a minimum $\bbeta_{1}(P_t)$.

According to Theorem~\ref{th:existence MM bounded rho} there exists at least one $\bbeta_1(P)$ that minimizes $R_P(\bbeta)$.
Now, suppose that $\bbeta_1(P)$ is unique.
Because $\bbeta_1(P)$ is regression equivariant, we may assume that $\bbeta_1(P)=\mathbf{0}$.
For the sake of brevity, let us write $\bbeta_{1,t}=\bbeta_1(P_t)$, $\bV_{0,t}=\bV_0(P_t)$, and~$R_t=R_{P_t}$.
From~\eqref{eq:conv Vt Yohai} it follows that for $t$ sufficiently large,
\begin{equation}
\label{eq:bounds Vt Yohai}
0<\lambda_k(\bV_0(P))/4\leq \lambda_k(\bV_{0,t})\leq\lambda_1(\bV_{0,t})\leq 4\lambda_1(\bV_0(P))<\infty.
\end{equation}
Now, consider a sequence $\{(\bbeta_{1,t},\bV_{0,t})\}$, such that $\|\bbeta_{1,t}\|\leq M$ and
$\bV_{0,t}$ satisfies~\eqref{eq:bounds Vt Yohai}.
Then the sequence $\{(\bbeta_{1,t},\bV_{0,t})\}$ lies in a compact set, so it has a convergent subsequence
$(\bbeta_{1,t_j},\bV_{0,t_j})\to(\bbeta_{1,L},\bV_0(P))$.
According to~\eqref{eq:prop Lemma 2 Yohai}, it follows that
\[
\begin{split}
\lim_{j\to\infty}
R_{t_j}(\bbeta_{1,t_j})
&=
\lim_{j\to\infty}
\int
\rho_1\left(
d(\bs,\bbeta_{1,t_j},\bV_{0,t_j})
\right)
\,\dd P_{t_j}(\bs)\\
&=
\int
\rho_1\left(
d(\bs,\bbeta_{1,L},\bV_0(P))
\right)
\,\dd P(\bs)
=
R_P(\bbeta_{1,L}).
\end{split}
\]
Now, suppose that $\bbeta_{1,L}\neq \mathbf{0}$.
Then, since $R_P(\bbeta)$ is uniquely minimized at $\bbeta=\mathbf{0}$,
this would mean that there exists $\eta>0$, such that together with~\eqref{eq:bound difference},
\[
R_{t_j}(\bbeta_{1,t_j})
>
R_P(\bbeta_{1,l})+2\eta
\geq
R_P(\mathbf{0})+2\eta
\geq
R_{t_j}(\mathbf{0})+\eta
>
R_{t_j}(\mathbf{0}),
\]
for $t_j$ sufficiently large,
This would mean that $\bbeta_{1,t_j}$ is not the minimizer of $R_{t_j}(\bbeta)$.
We conclude that $\bbeta_{1,L}=\mathbf{0}$, which proves the theorem.
\end{proof}
There are several examples of covariance functionals that satisfy~\eqref{eq:conv Vt Yohai},
such as the Minimum Covariance Determinant functional (see~\cite{cator-lopuhaa2012}) and the covariance S-functional
(see~\cite{lopuhaa1989}), including the Minimum Volume Ellipsoid functional.
For a structured covariance functional~$\bV(\btheta_0(P))$ to satisfy~\eqref{eq:conv Vt Yohai},
it is required that the mapping $\btheta\mapsto\bV(\btheta)$ is continuous.
This is true for all the examples mentioned in Section~\ref{sec:structured covariance model}.
In addition, the functional $\btheta_0(P)$ needs to be continuous.
An example is the S-functional $\btheta_0(P)$ defined in~\cite{lopuhaa-gares-ruizgazenARXIVE2022}.

A direct corollary of $\bbeta_1(P)$ being continuous,
is the consistency of the estimator~$\bbeta_{1,n}$.
\begin{corollary}
\label{cor:consistency MM-estimator bounded}
Suppose that $\rho_1:\R\to[0,\infty)$ satisfies (R-BND) and
suppose that $P$ is such that~(A) holds and that $\bX$ has full rank with probability one.
Suppose that $\bV_{0,n}\to\bV_0(P)$, with probability one.
Then for $n$ sufficiently large, there is at least one~$\bbeta_{1,n}$ that minimizes~$R_{n}(\bbeta)$,
with probability one.
If~$\bbeta_1(P)$ is the unique minimizer of $R_P(\bbeta)$, then for any sequence~$\bbeta_{1,n}$, $n=1,2\ldots$,
it holds that
\[
\lim_{n\to\infty}\bbeta_{1,n}=\bbeta_1(P),
\]
with probability one.
\end{corollary}
\begin{proof}
We apply Theorem~\ref{th:continuity MM bounded rho} to the sequence $\mathbb{P}_n$, $n=1,2,\ldots$,
of probability measures, where~$\mathbb{P}_n$ is the empirical measure corresponding to
$(\mathbf{y}_1,\mathbf{X}_1),\ldots,(\mathbf{y}_n,\mathbf{X}_n)$.
According to the Portmanteau Theorem (e.g., see Theorem~2.1 in~\cite{billingsley1968}),
$\mathbb{P}_n$ converges weakly to $P$, with probability one.
The corollary then follows from Theorem~\ref{th:continuity MM bounded rho}.
\end{proof}

For unbounded $\rho_1$, we cannot obtain continuity of the functional $\bbeta_1(P)$, for all sequences~$\{P_t\}$ that converge weakly to $P$.
However, we can establish strong consistency for the estimator~$\bbeta_{1,n}$.
\begin{theorem}
\label{th:consistency MM unbounded rho1}
Let $\rho_1:\R\to[0,\infty)$ satisfy (R-UNB).
Suppose that $\E_P\|\bs\|<\infty$ and that $\bX$ has full rank with probability one.
Suppose that $\bV_{0,n}\to\bV_0(P)$, with probability one,
and let~$\bbeta_{1,n}$ minimize~$R_{n}(\bbeta)$.
If $\bbeta_1(P)$ is the unique minimizer of $R_P(\bbeta)$, then
\[
\lim_{n\to\infty}\bbeta_{1,n}=\bbeta_1(P),
\]
with probability one.
\end{theorem}
\begin{proof}
For the sake of brevity, write $\bV_0$ instead of $\bV_0(P)$.
Since $\bV_{0,n}\to \bV_0$, with probability one, there exists
$0<L_1=\lambda_k(\bV_0)/4\leq 4\lambda_1(\bV_0)=L_2<\infty$,
such that, for $n$ sufficiently large, all eigenvalues of $\bV_{0,n}$ are between $L_1$ and $L_2$ with probability one.
Let $h(\bs;\bbeta,\bV)=\rho_1(\|\bV^{-1/2}(\by-\bX\bbeta)\|)$ and define
\[
\begin{split}
H(\bbeta,\bV)
&=
\int h(\bs;\bbeta,\bV)\,\dd P(\bs),\\
H_n(\bbeta,\bV)
&=
\int h(\bs;\bbeta,\bV)\,\dd \mathbb{P}_n(\bs).
\end{split}
\]
For $M>0$, consider the class of functions
$\mathcal{F}
=
\left\{
h(\cdot;\bbeta,\bV):
\|\bbeta\|\leq M;\, \lambda_k(\bV)\geq L_1
\right\}$.
Then, according to~\eqref{eq:prop rho general} and~\eqref{eq:bound rho1},
the class $\mathcal{F}$ has envelope
\[
\rho_1(s_0)+(\|\by\|+M\|\bX\|)a_1^{-1/2}\rho_1'(s_0),
\]
which is integrable, due to $\E_P\|\bs\|<\infty$.
Hence, by dominated convergence, $H(\bbeta,\bV)$ is continuous on the set
$K_M=\{(\bbeta,\bV): \|\bbeta\|\leq M;\, \lambda_k(\bV)\geq L_1\}$.

Moreover, the graphs of functions in $\mathcal{F}$ have polynomial discrimination.
This can be shown similar to the proof of Lemma~B.6 in~\cite{lopuhaa-gares-ruizgazenARXIVE2022}.
From Theorem~24 in~\cite{pollard1984}, we may then conclude
\begin{equation}
\label{eq:uniform strong law}
\sup_{(\bbeta,\bV)\in K_M}
\left|
H_n(\bbeta,\bV)-H(\bbeta,\bV)
\right|
\to0,
\end{equation}
with probability one.
As a first consequence, we find that
\begin{equation}
\label{eq:bound R-difference}
\left|
R_P(\mathbf{0})
-
R_n(\mathbf{0})
\right|
\leq
\left|
H_n(\mathbf{0},\bV_{0,n})
-
H(\mathbf{0},\bV_{0,n})
\right|
+
\left|
H(\mathbf{0},\bV_{0,n})
-
H(\mathbf{0},\bV_0)
\right|
\to 0,
\end{equation}
with probability one, due to~\eqref{eq:uniform strong law}
and continuity of~$H(\bbeta,\bV)$.
Next, we argue there exists~$M>0$, such that for $n$ sufficiently large $\|\bbeta_{1,n}\|\leq M$.
Since $\E_P\|\bs\|<\infty$, as in the proof of Theorem~\ref{th:existence MM unbounded rho}(i),
this ensures that
\[
R_P(\mathbf{0})
=
\int \rho_1\left(\|\bV_0^{-1/2}\by\|\right)\,\dd P(\bs)<\infty.
\]
Then, consider the set $A_M$ defined in~\eqref{def:AM}
and choose $M>0$, such that
\[
\rho_1(\sqrt{M}/(2\sqrt{\lambda_1}))P(A_M)>R_P(\textbf{0}).
\]
Then, for $n$ sufficiently large, we must have $\|\bbeta_{1,n}\|\leq M$,
since otherwise, according to~\eqref{eq:lower bound quadratic form},
\[
R_{n}(\bbeta_{1,n})
\geq
\rho_1(\sqrt{M}/(2\sqrt{\lambda_1}))\mathbb{P}_n(A_M)
\to
\rho_1(\sqrt{M}/(2\sqrt{\lambda_1}))P(A_M)
>
R_P(\mathbf{0}),
\]
as $n\to\infty$, with probability one,
which would imply that for $n$ sufficiently large,
$R_{n}(\bbeta_{1,n})>R_{n}(\mathbf{0})$, with probability one.

Then suppose that $\bbeta_1(P)$ is the unique minimizer of $R_P(\bbeta)$.
Because $\bbeta_1(P)$ is regression equivariant, we may assume that $\bbeta_1(P)=\mathbf{0}$.
This means that
for any $\delta>0$, there exist $\alpha>0$, such that
\[
\inf_{\|\bbeta\|>\delta}
\int
\rho_1\left(\frac{\|\bV_0^{-1/2}(\by-\bX\bbeta)\|}{1+\alpha}\right)\,\dd P(\bs)
>
R_P(\mathbf{0}).
\]
Because $\bV_{0,n}^{-1/2}\bV_0^{1/2}\to \bI_k$, with probability one,
we can choose $n$ sufficiently large such that $\lambda_k(\bV_{0,n}^{-1/2}\bV_0^{1/2})\geq 1/(1+\alpha)^2$.
Then, since
\[
\begin{split}
\|\bV_{0,n}^{-1/2}(\by-\bX\bbeta)\|^2
&=
\|\bV_{0,n}^{-1/2}\bV_0^{1/2}\bV_0^{-1/2}(\by-\bX\bbeta)\|^2\\
&\geq
\lambda_k(\bV_{0,n}^{-1/2}\bV_0^{1/2})
\|\bV_0^{-1/2}(\by-\bX\bbeta)\|^2
\geq
\frac{\|\bV_0^{-1/2}(\by-\bX\bbeta)\|^2}{(1+\alpha)^2},
\end{split}
\]
we find
\[
\inf_{\|\bbeta\|>\delta}
H_n(\bbeta,\bV_{0,n})
\geq
\inf_{\|\bbeta\|>\delta}
\int
\rho_1\left(\frac{\|\bV_0^{-1/2}(\by-\bX\bbeta)\|}{1+\alpha}\right)\,\dd P(\bs)
>
R_P(\mathbf{0}).
\]
Furthermore,
$\left|
R_P(\mathbf{0})
-
H_n(\mathbf{0},\bV_{0,n})
\right|
\to0$,
with probability one, as $n\to\infty$,
according to~\eqref{eq:bound R-difference}.
Hence, for $n$ sufficiently large, we would find that for all $\delta>0$,
\[
\inf_{\|\bbeta\|>\delta}
H_n(\bbeta,\bV_{0,n})
>
H_n(\mathbf{0},\bV_{0,n}).
\]
Therefore, for all $\delta>0$, we must have $\|\bbeta_{1,n}\|\leq \delta$, for $n$ sufficiently large,
with probability one.
This means $\bbeta_{1,n}\to\mathbf{0}$, with probability one.
\end{proof}
Asymptotic properties of the MM-estimator for linear mixed effects models in~\cite{copt&heritier2007}
was only considered for the simple model with a fixed design matrix and normal errors.
As a consequence of their Theorem~1, the MM-estimator would be weakly consistent.
However, the requirement that the covariance estimator used in the first step is consistent, seems to be missing.
Our Corollary~\ref{cor:consistency MM-estimator bounded} and Theorem~\ref{th:consistency MM unbounded rho1}
establish consistency for the MM-estimator for a larger class of linear mixed effects models
and under very mild conditions on the distribution~$P$.
Theorem~\ref{th:consistency MM unbounded rho1} is equivalent to Theorem~3.1 in~\cite{lopuhaa1992highly}
for multivariate location and scatter.
Corollary~\ref{cor:consistency MM-estimator bounded} extends this result to bounded $\rho$-functions.
Furthermore, Corollary~\ref{cor:consistency MM-estimator bounded} is obtained under conditions
on the distribution~$P$, that are much weaker than the ones for the MM-estimators considered
in~\cite{SalibianBarrera-VanAelst-Willems2006} and~\cite{kudraszow-maronna2011},
which restrict themselves to distributions with an elliptically contoured density.

\section{Global robustness: breakdown point}
\label{subsec:BDP Yohai}
Consider a collection of points $\mathcal{S}_n=\{\bs_i=(\by_i,\bX_i),i=1,\ldots,n\}\subset \R^k\times \R^{kq}$.
To emphasize the dependence on the collection $\mathcal{S}_n$,
we sometimes denote the estimators in Definition~\ref{def:MM-estimator general} by~$\bbeta_{1,n}(\mathcal{S}_n)$ and $\bV_{0,n}(\mathcal{S}_n)$.
To investigate the global robustness of $\bbeta_{1,n}$,
we compute that finite-sample (replacement) breakdown point.
For a given collection $\mathcal{S}_n$ the finite-sample breakdown point
(see~\cite{donoho&huber1983})
of regression estimator $\bbeta_{1,n}$ is defined as the smallest proportion of points
from~$\mathcal{S}_n$ that one needs to replace in order to
carry the estimator over all bounds.
More precisely,
\begin{equation}
\label{def:BDP beta}
\epsilon_n^*(\bbeta_{1,n},\mathcal{S}_n)
=
\min_{1\leq m\leq n}
\left\{
\frac{m}{n}:
\sup_{\mathcal{S}_m'}
\left\|
\bbeta_{1,n}(\mathcal{S}_n)-\bbeta_{1,n}(\mathcal{S}_m')
\right\|
=\infty
\right\},
\end{equation}
where the minimum runs over all possible collections $\mathcal{S}_m'$ that can be obtained from $\mathcal{S}_n$
by replacing~$m$ points of $\mathcal{S}_n$ by arbitrary points in $\R^k\times \R^{kq}$.
The finite sample (replacement) breakdown point of a covariance estimator $\bV_{0,n}$ at a collection~$\mathcal{S}_n$,
is defined as
\begin{equation}
\label{def:BDP V}
\epsilon_n^*(\bV_{0,n},\mathcal{S}_n)
=
\min_{1\leq m\leq n}
\left\{
\frac{m}{n}:
\sup_{\mathcal{S}_m'}
\text{dist}(\bV_{0,n}(\mathcal{S}_n),\bV_{0,n}(\mathcal{S}_m'))
=\infty
\right\},
\end{equation}
with $\text{dist}(\cdot,\cdot)$ defined as
$\text{dist}(\bA,\mathbf{B})
=
\max\left\{
\left|\lambda_1(\bA)-\lambda_1(\mathbf{B})\right|,
\left|\lambda_k(\bA)^{-1}-\lambda_k(\mathbf{B})^{-1}\right|
\right\}$,
where the minimum runs over all possible collections $\mathcal{S}_m'$ that can be obtained from $\mathcal{S}_n$
by replacing~$m$ points of $\mathcal{S}_n$ by arbitrary points in $\R^k\times \R^{kq}$.
So the breakdown point of $\bV_{0,n}$ is the smallest proportion of points from~$\mathcal{S}_n$ that one needs to replace in order to
make the largest eigenvalue of~$\bV_{0,n}(\mathcal{S}_m')$ arbitrarily large (explosion), or
to make the smallest eigenvalue of~$\bV_{0,n}(\mathcal{S}_m')$ arbitrarily small (implosion).
When we estimate a structured covariance matrix $\bV(\btheta)$ by $\bV(\btheta_{0,n})$,
we need to specify what the breakdown point of $\btheta_{0,n}$ is.
Since, the estimator $\btheta_{0,n}$ determines the covariance estimator $\bV(\btheta_{0,n})$,
it seems natural to let the breakdown point of $\btheta_{0,n}$ correspond to the breakdown point of the covariance estimator
(see also~\cite{lopuhaa-gares-ruizgazenARXIVE2022}).

The breakdown behavior of $\bbeta_{1,n}$ depends on whether the function $\rho_1$ in Definition~\ref{def:MM-estimator general}
is bounded or unbounded.
\cite{huber1984} pointed put that the breakdown point of location M-estimators constructed with a bounded $\rho$-function
not only depends on the function $\rho$, but also on the configuration of the sample.
Depending on the configuration of the sample, the breakdown point can be any value between 0 and $1/2$.
For this reason, the multivariate location M-estimator (see~\cite{lopuhaa1992highly}) is constructed with an unbounded $\rho$-function.
In this way, the location M-estimator inherits the breakdown point of the initial covariance estimator.

Unfortunately, the use of an unbounded function $\rho_1$ in Definition~\ref{def:MM-estimator general},
does not seem suitable for the breakdown behavior of the regression MM-estimator.
The presence of the design matrices~$\bX_i$ makes things more complicated than in the multivariate location case.
Nevertheless, for unbounded~$\rho_1$, we can establish a result similar to the one in~\cite{lopuhaa1992highly},
when all design matrices are the same.
An example is of course when all $\bX_i=\bI_k$ as in the location-scale model,
but another example occurs in linear mixed effects models for which
all subjects have the same the design matrix $\bX$ representing particular contrasts for the fixed effects.
\begin{proposition}
\label{prop:BDP MM unbounded}
Suppose that $\rho_1$ satisfies (R-UNB).
Let $\mathcal{S}_n\subset \R^{p}$ be a collection of~$n$ points $\bs_i=(\by_i,\bX_i)$,  $i=1,\ldots,n$.
Suppose that $\bX_i=\bX$, for all $i=1,\ldots,n$, where $\bX$ is fixed and has full rank.
Then for any $\bbeta_{1,n}$ that minimizes $R_n(\bbeta)$,
it holds that
\[
\epsilon_n^*(\bbeta_{1,n},\mathcal{S}_n)
\geq
\epsilon_n^*(\bV_{0,n},\mathcal{S}_n).
\]
\end{proposition}
\begin{proof}
For $\bt\in\R^k$, let
\[
\widetilde{R}_n(\bt)
=
\frac1n
\sum_{i=1}^{n}
\left\{
\rho_1\left(
\sqrt{(\by_i-\bt)^T\bV_{0,n}^{-1}(\by_i-\bt)}
\right)
-
\rho_1\left(
\sqrt{\by_i^T\bV_{0,n}^{-1}\by_i}
\right)
\right\},
\]
be the object function for the location M-estimator in~\cite{lopuhaa1992highly}.
Then we can write
\[
R_n(\bbeta)
=
R_n(\mathbf{0})
+
\widetilde{R}_n(\bX\bbeta).
\]
Because $\bX$ has full rank, $\bbeta_{1,n}$ minimizes $R_n(\bb)$ if and only if
$\bt_{1,n}=(\bX^T\bX)^{-1}\bX^T\bbeta_{1,n}$ minimizes~$\widetilde{R}_n(\bt)$.
As $\bX$ is considered to be fixed, this means that $\bbeta_{1,n}$ breaks down precisely when~$\bt_{1,n}$ does.
Hence from Theorem~4.1 in~\cite{lopuhaa1992highly} we conclude that
$\epsilon_n^*(\bbeta_{1,n},\mathcal{S}_n)\geq\epsilon_n^*(\bV_{0,n},\mathcal{S}_n)$.
\end{proof}
The use of a bounded function $\rho_1$ in Definition~\ref{def:MM-estimator general} also does not seem very suitable,
in view of the results found by~\cite{huber1984}.
However, the approach followed by~\cite{yohai1987}, which relates the bounded function $\rho_1$
to another bounded function $\rho_0$ in the first stage,
turns out to be adequate.
Let $\bbeta_{0,n}$ be an initial regression estimate, such that together with the covariance estimate $\bV_{0,n}$,
it holds that
\begin{equation}
\label{eq:cond stage 1}
\frac1n
\sum_{i=1}^n
\rho_0
\left(
\sqrt{(\by_i-\bX_i\bbeta_{0,n})^T\bV_{0,n}^{-1}(\by_i-\bX_i\bbeta_{0,n})}
\right)
=
b_0,
\end{equation}
for a function $\rho_0$ that satisfies (R-BND)
and suppose that $\rho_1$ satisfies (R-BND), such that
\begin{equation}
\label{eq:ineq rho}
\frac{\rho_1(s)}{a_1}
\leq
\frac{\rho_0(s)}{a_0}.
\end{equation}
Next, we proceed as in Definition~\ref{def:MM-estimator general}, i.e., define $\bbeta_{1,n}$ as
the vector that minimizes~\eqref{eq:MM estimator general}.
Estimates~$\bbeta_{0,n}$ and $\bV_{0,n}$ can be any two (high breakdown) estimators satisfying~\eqref{eq:cond stage 1}.
However, natural candidates for our setup are the S-estimates $(\bbeta_{0,n},\btheta_{0,n})$, with $\bV_{0,n}=\bV(\btheta_{0,n})$,
defined in~\cite{lopuhaa-gares-ruizgazenARXIVE2022} by means of the function $\rho_0$.

In order to formulate the breakdown point of $\bbeta_{1,n}$ using bounded $\rho$-functions, we first need to discuss the following.
Recall that $(\by_1,\bX_1),\ldots,(\by_n,\bX_n)$ are represented as points in $\R^k\times\R^{kq}$.
Note however, that for linear models with intercept the first column of each $\bX_i$ consists of 1's.
This means that the points $(\by_i,\bX_i)$ are concentrated in a lower dimensional subset of $\R^k\times\R^{kq}$.
A similar situation occurs when all $\bX_i$ are equal to the same design matrix.
In view of this, define~$\mathcal{X}\subset\R^{kq}$ as the subset with the lowest dimension
$p=\text{dim}(\mathcal{X})\leq kq$ satisfying
$P(\bX\in \mathcal{X})=1$.
Hence, $P$ is concentrated on the subset $\R^k\times \mathcal{X}$ of $\R^k\times\R^{kq}$, which
may be of lower dimension~$k+p$ than $k+kq$.
Let $\mathcal{S}_n=\{\bs_1,\ldots,\bs_n\}$, with $\bs_i=(\by_i,\bX_i)$ be a collection of $n$ points
in~$\R^k\times \mathcal{X}$.
Define
\begin{equation}
\label{def:k(S)}
\kappa(\mathcal{S}_n)
=
\text{maximal number of points of $\mathcal{S}_n$ lying on the same hyperplane in~$\R^k\times \mathcal{X}$.}
\end{equation}
For example, if the distribution $P$ is absolutely continuous, then
$\kappa(\mathcal{S}_n)\leq k+p$ with probability one.
We then have the following theorem.
\begin{theorem}
\label{th:BDP MM Yohai}
Suppose that $(\bbeta_{0,n},\bV_{0,n})$ satisfy~\eqref{eq:cond stage 1}, with $\rho_0$ such that~\eqref{eq:ineq rho} holds,
where~$\rho_0$ and~$\rho_1$ satisfy (R-BND).
Let $\mathcal{S}_n\subset \R^{k}\times \mathcal{X}$ be a collection of~$n$ points $\bs_i=(\by_i,\bX_i)$,  $i=1,\ldots,n$.
Let $r_0=b_0/a_0$ and suppose that $0<r_0\leq (n-k(\mathcal{S}_n))/(2n)$,
where $k(\mathcal{S}_n)$ is defined by~\eqref{def:k(S)}.
Then for any $\bbeta_{1,n}$ that minimizes $R_n(\bbeta)$,
it holds that
\[
\epsilon_n^*(\bbeta_{1,n},\mathcal{S}_n)
\geq
\min\left(
\epsilon_n^*(\bV_{0,n},\mathcal{S}_n),\frac{\lceil{nr_0}\rceil}n
\right).
\]
\end{theorem}
\begin{proof}
Suppose we replace $m$ points, where $m$ is such that
\[
m\leq \lceil nr_0\rceil-1
\quad\text{and}\quad
m\leq n\epsilon_n^*(\bV_{0,n},\mathcal{S}_n)-1.
\]
Let $\mathcal{S}_m'$ be the corrupted collection of points.
Write $\bbeta_{0,m}=\bbeta_{0,n}(\mathcal{S}_m')$, and $\bV_{0,m}=\bV_{0,n}(\mathcal{S}_m')$.
Then $\bV_{0,m}$ does not break down, so that there exist constants $0<L_1\leq L_2<\infty$,
not depending on $\mathcal{S}_m'$ such that
$0<L_1\leq\lambda_k(\bV_{0,m})\leq \lambda_1(\bV_{0,m})\leq L_2<\infty$.
For any $\bbeta\in\R^q$,
define the cylinder
$\mathcal{C}_{1,m}(\bbeta)
=
\left\{
(\by,\bX)\in\R^{p}:
(\by-\bX\bbeta)^T
\bV_{0,m}^{-1}
(\by-\bX\bbeta)
\leq
c_1^2
\right\}$.
Consider the function
\[
R_m(\bbeta)
=
\frac{1}{n}
\sum_{\bs_i\in S_m'}
\rho_1\left(
\sqrt{(\by_i-\bX_i\bbeta)^T\bV_{0,m}^{-1}(\by_i-\bX_i\bbeta)}
\right),
\]
for the corrupted sample $\mathcal{S}_m'$.
For any $\bbeta$ that minimizes $R_m(\bbeta)$, it holds $R_m(\bbeta)\leq R_m(\bbeta_{0,m})$.
Therefore, for such $\bbeta$,
according to~\eqref{eq:ineq rho} and~\eqref{eq:cond stage 1}, we have that
\[
\begin{split}
R_m(\bbeta)
&\leq
\frac{1}{n}
\sum_{\bs_i\in S_m'}
\rho_1\left(
\sqrt{(\by_i-\bX_i\bbeta_{0,m})^T\bV_{0,m}^{-1}(\by_i-\bX_i\bbeta_{0,m})}
\right)\\
&\leq
\frac{a_1}{a_0}
\frac{1}{n}
\sum_{\bs_i\in S_m'}
\rho_0\left(
\sqrt{(\by_i-\bX_i\bbeta_{0,m})^T\bV_{0,m}^{-1}(\by_i-\bX_i\bbeta_{0,m})}
\right)
=
r_0a_1.
\end{split}
\]
Let $\mathbb{P}_m'$ be the empirical measure corresponding to the corrupted collection $\mathcal{S}_m'$.
Then it holds that
\[
\begin{split}
\mathbb{P}_m'(\mathcal{C}_{1,m}(\bbeta))
&=
\frac1n
\sum_{\bs_i\in S_m'}
\mathds{1}
\left\{
\bs_i\in
\mathcal{C}_{1,m}(\bbeta)
\right\}\\
&\geq
1-
\frac{1}{na_1}
\sum_{\bs_i\in S_m'}
\rho_1\left(
\sqrt{(\by_i-\bX_i\bbeta)^T\bV_{0,m}^{-1}(\by_i-\bX_i\bbeta)}
\right)
\geq
1-r_0.
\end{split}
\]
It follows that the cylinder
$\mathcal{C}_{1,m}(\bbeta)$
must contain at least $\lceil n-nr_0\rceil$ number of points from the corrupted collection
$\mathcal{S}'_m$.
Furthermore, since $r_0\leq (n-k(\mathcal{S}_n))/(2n)$,
for any such subset of $\mathcal{S}'_m$ it holds that
it contains
$\lceil n-nr_0\rceil-m
=
n-\lfloor{nr_0}\rfloor-\lceil nr_0\rceil+1
\geq
k(\mathcal{S}_n)+1$
points of the original collection $\mathcal{S}_n$.
Let $J_0$ be a subset of
$k(\mathcal{S}_n)+1$ points from the original collection $\mathcal{S}_n$
contained in $\mathcal{C}_{1,m}(\bbeta)$.
By definition, $k(\mathcal{S}_n)+1$ original points cannot be on the same hyperplane, so that
\[
\gamma_n
=
\inf_{J\subset \mathcal{S}_n}
\inf_{\|\bgamma\|=1}\max_{\bs\in J}
\|\bX\bgamma\|>0.
\]
where the first infimum runs over all subsets $J\subset \mathcal{S}_n$ of $k(\mathcal{S}_n)+1$ points.
By definition of~$\gamma_n$,
there exists an original point~$\bs_0\in J_0\subset \mathcal{S}_n\cap \mathcal{C}_{1,m}(\bbeta)$, such that
\[
\|\bbeta\|
=
\|\bX_0\bbeta\|
\times
\frac{\|\bbeta\|}{\|\bX_0\bbeta\|}
\leq
\frac{1}{\gamma_n}
\|\bX_0\bbeta\|.
\]
Because $\bs_0\in \mathcal{C}_{1,m}(\bbeta)$, it follows
that
\[
\|\by_0-\bX_0\bbeta\|^2
\leq
(\by_0-\bX_0\bbeta)^T\bV_{0,n}^{-1}(\by_0-\bX_0\bbeta)
\lambda_1(\bV_{0,m})
\leq
c_1^2L_2,
\]
and because $\bs_0\in \mathcal{S}_n$, we have that
\[
\|\bX_0\bbeta\|
\leq
c_1\sqrt{L_2}
+
\max_{(\by_i,\bX_i)\in \mathcal{S}_n}
\|\by_i\|
<\infty.
\]
We conclude that for minimizing $R_m(\bbeta)$ we can restrict ourselves to a compact set $K_n$, only depending on the original collection~$\mathcal{S}_n$.
Firstly, since $R_m(\bbeta)$ is continuous, this implies there exists at least one $\bbeta_{1,n}(\mathcal{S}'_m)$,
which minimizes $R_m(\bbeta)$.
Secondly, since any $\bbeta_{1,n}(\mathcal{S}'_m)$ must be in~$K_n$, which only depends on the original collection~$\mathcal{S}_n$,
the estimate $\bbeta_{1,n}(\mathcal{S}'_m)$ does not break down.
\end{proof}
Theorem~\ref{th:BDP MM Yohai} is comparable to Theorem~1 in~\cite{SalibianBarrera-VanAelst-Willems2006} and
Theorem~3 in~\cite{kudraszow-maronna2011}
for MM-estimators for multivariate location and scatter and for multivariate linear regression, respectively.
The breakdown point for MM-estimators for linear mixed effects models has only been discussed in~\cite{copt&heritier2007}.
They conjecture that the exact value can be derived using the technique in~\cite{vanaelst&willems2005},
but do not pursue a rigorous derivation.
Together with Proposition~\ref{prop:BDP MM unbounded}, the result in Theorem~\ref{th:BDP MM Yohai} provides
sufficient conditions for the MM-estimators in the linear mixed effects model used in~\cite{copt&heritier2007},
to inherit the breakdown point from the initial covariance estimate.

It can be shown that if $\bV_{0,n}$ satisfies~\eqref{eq:cond stage 1} for some $\bbeta_{0,n}$,
it must have a breakdown point that is less than or equal to $\lceil nr_0\rceil/n$
(e.g., see the proof of Theorem 4 in~\cite{lopuhaa-gares-ruizgazenARXIVE2022}).
This means that from Theorem~\ref{th:BDP MM Yohai}, we have
$\epsilon_n^*(\bbeta_{1,n},\mathcal{S}_n)
\geq
\epsilon_n^*(\bV_{0,n})$.
Moreover, if $(\bbeta_{0,n},\btheta_{0,n})$ are S-estimators, as defined in~\cite{lopuhaa-gares-ruizgazenARXIVE2022},
such that $\bV_{0,n}=\bV(\btheta_{0,n})$,
then $\epsilon_n^*(\bV_{0,n})=\lceil nr_0\rceil/n$ according to
Theorem~4 in~\cite{lopuhaa-gares-ruizgazenARXIVE2022}, so that
\[
\epsilon_n^*(\bbeta_{1,n},\mathcal{S}_n)
\geq
\frac{\lceil{nr_0}\rceil}n.
\]
The largest possible value of the breakdown point occurs when $r_0=(n-\kappa(\mathcal{S}_n))/(2n)$, in which case
$\lceil nr_0\rceil/n=\lceil (n-\kappa(\mathcal{S}_n))/2\rceil/n=\lfloor (n-\kappa(\mathcal{S}_n)+1)/2\rfloor/n$.
When the collection~$\mathcal{S}_n$ is in general position, then $\kappa(\mathcal{S}_n)= k+p$.
In that case the breakdown point is at least~$\lfloor (n-k-p+1)/2\rfloor/n$.
When all $\bX_i$ are equal to the same $\bX$, one has $p=0$ and $\kappa(\mathcal{S}_n)=k$.
In that case, the breakdown point is at least $\lfloor (n-k+1)/2\rfloor/n$.
This coincides with the maximal breakdown point for affine equivariant estimators
for $k\times k$ covariance matrices (see~\cite{davies1987}).

\section{Score equations}
\label{sec:score equations}
Recall the definition of the functional $\bbeta_1(P)$ in Section~\ref{subsec:def MM Yohai},
which minimizes $R_{P}(\bbeta)$.
Then $\bbeta_1(P)$ is also a solution of $\partial R_P(\bbeta)/\partial \bbeta=\mathbf{0}$.
In order to allow changing the order of integration and differentiation in $R_P(\bbeta)$,
we require an additional condition on $\rho_1$.
\begin{quote}
\begin{itemize}
\item[(R-CD1)]
$\rho_1$ is continuously differentiable and $u_1(s)=\rho_1'(s)/s$ is continuous,
\end{itemize}
\end{quote}
If $\rho_1$ satisfies (R-CD1), then
\[
\frac{\partial}{\partial \bbeta}
\rho_1\left(
\sqrt{(\by-\bX\bbeta)^T\bV_0(P)^{-1}(\by-\bX\bbeta)}
\right)
=
u_1(d_0)\bX^T\bV_0(P)^{-1}(\by-\bX\bbeta),
\]
where $d_0=d(\bs,\bbeta,\bV_0(P))$, as defined in~\eqref{def:mahalanobis}.
This means that
\[
\left\|
\frac{\partial}{\partial \bbeta}
\rho_1\left(
\sqrt{(\by-\bX\bbeta)^T\bV_0(P)^{-1}(\by-\bX\bbeta)}
\right)
\right\|
=
\frac{|u_1(d_0)d_0|\cdot\|\bX\|}{\sqrt{\lambda_k(\bV_0(P))}}.
\]
When $\rho_1$ satisfies either (R-BND) or (R-UNB),
the function $u_1(s)s$ is uniformly bounded.
This means that in both cases the right hand side is bounded by a constant times $\|\bX\|$.
Hence, if~$\E_P\|\bX\|<\infty$, then by dominated convergence
\[
\begin{split}
&
\frac{\partial}{\partial \bbeta}
\int
\rho_1\left(
\sqrt{(\by-\bX\bbeta)^T\bV_0(P)^{-1}(\by-\bX\bbeta)}
\right)
\,\dd
P(\by,\bX)\\
&=
\int
\frac{\partial}{\partial \bbeta}
\rho_1\left(
\sqrt{(\by-\bX\bbeta)^T\bV_0(P)^{-1}(\by-\bX\bbeta)}
\right)
\,\dd
P(\by,\bX)\\
&=
\int
u_1(d_0)\bX^T\bV_0(P)^{-1}(\by-\bX\bbeta)
\,\dd
P(\by,\bX).
\end{split}
\]
We conclude that, if $\E_P\|\bX\|<\infty$, the functional $\bbeta_1(P)$ satisfies
score equations
\begin{equation}
\label{eq:score psi MM Yohai}
\int
\Psi_1(\bs,\bbeta,\bV_0(P))
\,\dd
P(\by,\bX)
=
\mathbf{0},
\end{equation}
where
\begin{equation}
\label{def:psi beta MM general}
\Psi_1(\bs,\bbeta,\bV)
=
u_1(d)\bX^T\bV^{-1}(\by-\bX\bbeta),
\end{equation}
with $d=d(\bs,\bbeta,\bV)$, as defined in~\eqref{def:mahalanobis}.

Score equation~\eqref{eq:score psi MM Yohai} coincides with equation (3.8) in~\cite{lopuhaa1992highly}
for the multivariate location-scale model.
If $P$ is the empirical measure $\mathbb{P}_n$ corresponding to $(\by_1,\bX_1),\ldots,(\by_n,\bX_n)$,
then~\eqref{eq:score psi MM Yohai} coincides with equation (2.6) for the multiple regression model in~\cite{yohai1987}
and with equation~(8) for the linear mixed effects model~\eqref{def:linear mixed effects model Copt}
in~\cite{copt&heritier2007}.
Furthermore, for the empirical measure~$\mathbb{P}_n$,
equation~\eqref{eq:score psi MM Yohai} is also similar to equation~(16) for the location MM-estimator
in~\cite{SalibianBarrera-VanAelst-Willems2006}
and to equation~(2.10) for the multivariate linear regression MM-estimator in~\cite{kudraszow-maronna2011}.

Let
\begin{equation}
\label{def:Lambda beta MM general}
\Lambda_1(\bbeta,\bV)
=
\int \Psi_1(\bs,\bbeta,\bV)\,\dd P(\bs).
\end{equation}
A vector $\bb\in\R^q$ is called a \emph{point of symmetry} of $P$, if for almost all $\bX$, it holds that
\[
P\left(\bX\bb+A\mid\bX\right)
=
P\left(\bX\bb-A\mid\bX\right),
\]
for all measurable sets $A\subset\R^k$,
where for $\lambda\in\R$ and $\bb\in\R^q$,
$\bX\bb+\lambda A$ denotes the set~$\{\bX+\lambda \by:\by\in A\}$.
If $\bb$ is a point of symmetry of $P$, it has the property that
\begin{equation}
\label{eq:point of symmetry}
\Lambda_1(\bb,\bV)=\mathbf{0}
\quad
\text{for all non-singular symmetric matrices }\bV.
\end{equation}
This will become very useful in determining asymptotic properties, such as the influence function for $\bbeta_1(P)$ and
asymptotic normality of $\bbeta_{1,n}$.
Note that if $P$ is such that $\by\mid\bX$ has an elliptically contoured density as defined in~\eqref{eq:elliptical}, then
the vector $\bmu=\bX\bbeta$ is a point of symmetry.

\section{Local robustness: the influence function}
\label{subsec:IF MM Yohai}
For $0<\epsilon<1$ and $\bs_0=(\by_0,\bX_0)\in\R^{k}\times \R^{kq}$ fixed, consider
$P_{\epsilon,\bs_0}=(1-\epsilon)P+\epsilon\delta_{\bs_0}$, as defined in~\eqref{def:perturbed P}.
The \emph{influence function} of the functional~$\bbeta_1(\cdot)$ at probability measure $P$,
is defined as
\begin{equation}
\label{def:IF regression Yohai}
\text{IF}(\bs;\bbeta_1,P)
=
\lim_{\epsilon\downarrow0}
\frac{\bbeta_1((1-\epsilon)P+\epsilon\delta_{\bs_0})-\bbeta_1(P)}{\epsilon},
\end{equation}
if this limit exists (see~\cite{hampel1974}).

We intend to include both bounded and unbounded functions $\rho_1$ in Definition~\ref{def:MM-functional general}.
For bounded~$\rho_1$, it follows from Theorem~\ref{th:continuity MM bounded rho} that,
under suitable conditions,
the functional~$\bbeta_1(P)$ is continuous.
In particular, this means that
\begin{equation}
\label{eq:cont beta1}
\lim_{\epsilon\downarrow0}
\bbeta_1(P_{\epsilon,\bs_0})=\bbeta_1(P).
\end{equation}
For unbounded $\rho_1$, the functional $\bbeta_1(P)$ is not necessarily continuous, but we can still establish~\eqref{eq:cont beta1}.
\begin{lemma}
\label{lem:beta1 continuous}
Let $\rho_1:\R\to[0,\infty)$ satisfy (R-UNB).
Suppose that $\E_P\|\bs\|<\infty$ and that $\bX$ has full rank with probability one.
Suppose that $\bV_0(P_{\epsilon,\bs_0})$ exists and that $\bV_0(P_{\epsilon,\bs_0})\to\bV_0(P)$, as $\epsilon\downarrow0$.
Suppose that~$\bbeta_{1}(P_{\epsilon,\bs_0})$ minimizes $R_{P_{\epsilon,\bs_0}}(\bbeta)$.
If $\bbeta_1(P)$ is the unique minimizer of~$R_P(\bbeta)$, then
\[
\lim_{\epsilon\downarrow0}
\bbeta_1(P_{\epsilon,\bs_0})=\bbeta_1(P).
\]
\end{lemma}
\begin{proof}
Define the functions $h(\bs;\bbeta,\bV)$ and $H(\bbeta,\bV)$ as in the proof of Theorem~\ref{th:consistency MM unbounded rho1},
and let
\[
H_{\epsilon,\bs_0}(\bbeta,\bV)=\int h(\bs;\bbeta,\bV)\,\dd P_{\epsilon,\bs_0}(\bs).
\]
Let $K_M$ be the same set of pairs $(\bbeta,\bV)$ as in the proof of Theorem~\ref{th:consistency MM unbounded rho1}.
Then for $(\bbeta,\bV)\in K_M$,
\[
|h(\bs;\bbeta,\bV)|
\leq
\rho_1(s_0)+(\|\by\|+M\|\bX\|)L_1^{-1/2}\rho_1'(s_0).
\]
Instead of~\eqref{eq:uniform strong law}, we now have
\[
\sup_{(\bbeta,\bV)\in K_M}
\left|
H_{\epsilon,\bs_0}(\bbeta,\bV)-H(\bbeta,\bV)
\right|
\leq
\epsilon
\left(
\sup_{(\bbeta,\bV)\in K_M}
|H(\bbeta,\bV)|
+
\sup_{(\bbeta,\bV)\in K_M}
h(\bs_0;\bbeta,\bV)
\right).
\]
Because $\E_P\|\bs\|<\infty$, the first supremum on the right hand side is bounded
and similarly, the second supremum is bounded by a constant depending on $\bs_0$ and $\rho_1$.
Therefore,\[
\sup_{(\bbeta,\bV)\in K_M}
\left|
H_{\epsilon,\bs_0}(\bbeta,\bV)-H(\bbeta,\bV)
\right|
\to
0,\quad
\text{as }\epsilon\downarrow0.
\]
From here on, one can mimic the proof of Theorem~\ref{th:consistency MM unbounded rho1}
and show that $\bbeta_1(P_{\epsilon,\bs_0})\to\bbeta_1(P)$.
\end{proof}
Now, that we have established~\eqref{eq:cont beta1} for both bounded and unbounded $\rho_1$,
we have the following general result for the influence function.

\begin{theorem}
\label{th:IF MM}
Suppose that $\rho_1$ either satisfies (R-BND) and (R-CD1) or satisfies~(R-UNB), and
suppose that $\E_P\|\bs\|<\infty$.
Let $\bbeta_1(P_{\epsilon,\bs_0})$ and $\bbeta_1(P)$ be a solution to the minimization problem
in Definition~\ref{def:MM-functional general} at $P_{\epsilon,\bs_0}$ and $P$, respectively.
Suppose that
$(\bbeta_1(P_{\epsilon,\bs_0}),\bV_0(P_{\epsilon,\bs_0}))\to (\bbeta_1(P),\bV_0(P))$,
as $\epsilon\downarrow0$,
and suppose that $\bbeta_1(P)$ is a point of symmetry of $P$.
Suppose that~$\Lambda_1$, as defined in~\eqref{def:Lambda beta MM general},
has a partial derivative $\partial\Lambda_1/\partial\bbeta$ that is continuous
at $(\bbeta_1(P),\bV_0(P))$ and that
$\mathbf{D}_1=(\partial\Lambda_1/\partial\bbeta)(\bbeta_1(P),\bV_0(P))$ is non-singular.
Then for $\bs_0\in\R^p$,
\[
\mathrm{IF}(\bs_0;\bbeta_1,P)
=
-\mathbf{D}_1^{-1}\Psi_1(\bs_0,\bbeta_1(P),\bV_0(P)),
\]
where $\Psi_1$ is defined in~\eqref{def:psi beta MM general}.
\end{theorem}
\begin{proof}
Denote
\[
\begin{split}
\bxi_{\epsilon,\bs_0}
&=
(\bbeta_{1,\epsilon,\bs_0},\bV_{0,\epsilon,\bs_0})=(\bbeta_1(P_{\epsilon,\bs_0}),\bV_0(P_{\epsilon,\bs_0}))\\
\bxi_P
&=
(\bbeta_{1,P},\bV_{0,P})=(\bbeta_1(P),\bV_0(P))
\end{split}
\]
Since $\bV_{0,\epsilon,\bs_0}\to \bV_{0,P}$ and the fact that $\E_P\|\bs\|<\infty$,
it follows from Theorem~\ref{th:continuity MM bounded rho} and Lemma~\ref{lem:beta1 continuous}
that there exists $\bbeta_1(P_{\epsilon,\bs_0})$ minimizing $R_P(\bbeta)$ at $P=P_{\epsilon,\bs_0}$.
According to Section~\ref{sec:score equations}, this means that~$\bxi_{\epsilon,\bs_0}$ satisfies the score equation~\eqref{eq:score psi MM Yohai} for the
regression M-functional $\bbeta_1$ at $P_{\epsilon,\bs_0}$, that is
\[
\int \Psi_1(\bs,\bxi_{\epsilon,\bs_0})\,\dd P_{\epsilon,\bs_0}(\bs)=\mathbf{0}.
\]
We decompose as follows
\[
\begin{split}
\mathbf{0}
&=
\int \Psi_1(\bs,\bxi_{\epsilon,\bs_0})\,\dd P_{\epsilon,\bs_0}(\bs)\\
&=
(1-\epsilon)\int \Psi_1(\bs,\bxi_{\epsilon,\bs_0})\,\dd P(\bs)+\epsilon\Psi_1(\bs_0,\bxi_{\epsilon,\bs_0})
=
(1-\epsilon)\Lambda_1(\bxi_{\epsilon,\bs_0})+\epsilon\Psi_1(\bs_0,\bxi_{\epsilon,\bs_0}).
\end{split}
\]
We first determine the order of $\bbeta_{1,\epsilon,\bs_0}-\bbeta_{1,P}$, as $\epsilon\downarrow0$.
Because $\bxi\mapsto\Psi_1(\bs_0,\bxi)$ is continuous, it follows that
$\Psi_1(\bs_0,\bxi_{\epsilon,\bs_0})=\Psi_1(\bs_0,\bxi_P)+o(1)$,
as $\epsilon\downarrow0$.
Furthermore, because $\Lambda_1$ has a partial derivative $\partial\Lambda_1/\partial\bbeta$ that is continuous
at $\bxi_P=(\bbeta_{1,P},\bV_{0,P})$,
we have that
\[
\begin{split}
\Lambda_1(\bxi_{\epsilon,\bs_0})
&=
\Lambda_1(\bbeta_{1,P},\bV_{0,\epsilon,\bs_0})
+
\frac{\partial \Lambda_1}{\partial \bbeta}
(\bbeta_{1,P},\bV_{0,\epsilon,\bs_0})
(\bbeta_{1,\epsilon,\bs_0}-\bbeta_{1,P})
+
o(\|\bbeta_{1,\epsilon,\bs_0}-\bbeta_{1,P}\|)\\
&=
\Lambda_1(\bbeta_{1,P},\bV_{0,\epsilon,\bs_0})
+
\Big(
\bD_1+o(1)
\Big)
\Big(\bbeta_{1,\epsilon,\bs_0}-\bbeta_{1,P}\Big)
+
o(\|\bbeta_{1,\epsilon,\bs_0}-\bbeta_{1,P}\|).
\end{split}
\]
Since $\bbeta_{1,P}$ is a point of symmetry of $P$, according to~\eqref{eq:point of symmetry} it holds that
$\Lambda_1(\bbeta_{1,P},\bV_{0,\epsilon,\bs_0})=\mathbf{0}$.
It follows that
\[
\mathbf{0}=
(1-\epsilon)
\bD_1
\Big(\bbeta_{1,\epsilon,\bs_0}-\bbeta_{1,P}\Big)
+
o(\|\bbeta_{1,\epsilon,\bs_0}-\bbeta_{1,P}\|)
+
o(\epsilon)
+
\epsilon\Psi_1(\bs_0,\bxi_{P}).
\]
Because $\bD_1$ is non-singular
and $\Psi_1(\bs_0,\bxi_P)$ is fixed, this implies
$\bbeta_{1,\epsilon,\bs_0}-\bbeta_{1,P}=O(\epsilon)$.
After inserting this in the previous equality, it follows that
\[
\begin{split}
\mathbf{0}
&=
(1-\epsilon)
\bD_1
\Big(\bbeta_{1,\epsilon,\bs_0}-\bbeta_{1,P}\Big)
+
\epsilon\Psi_1(\bs_0,\bxi_P)
+o(\epsilon)\\
&=
\bD_1
\Big(\bbeta_{1,\epsilon,\bs_0}-\bbeta_{1,P}\Big)
+
\epsilon\Psi_1(\bs_0,\bxi_P)
+o(\epsilon).
\end{split}
\]
We conclude
\[
\frac{\bbeta_{1,\epsilon,\bs_0}-\bbeta_{1,P}}{\epsilon}
=
-\bD_1^{-1}\Psi_1(\bs_0,\bxi_P)+o(1).
\]
This means that the limit of the left hand side exists and
\[
\text{IF}(\bs_0;\bbeta_1,P)
=
\lim_{\epsilon\downarrow0}
\frac{\bbeta_1((1-\epsilon)P+\epsilon\delta_{\bs_0})-\bbeta_1(P)}{\epsilon}
=
-\bD_1^{-1}\Psi_1(\bs_0,\bxi_P).
\]
\end{proof}
When $P$ is such that $\by\mid\bX$ has an elliptically contoured density~\eqref{eq:elliptical}
we can obtain a more detailed expression for the influence function.
This requires the following additional condition on the function~$\rho_1$.
\begin{quote}
\begin{itemize}
\item[(R-CD2)]
$\rho_1$ is twice continuously differentiable.
\end{itemize}
\end{quote}
We then have the following corollary.

\begin{corollary}
\label{cor:IF MM}
Suppose that $P$ is such that $\by\mid\bX$ has an elliptically contoured
density~$f_{\bmu,\bSigma}$ from~\eqref{eq:elliptical},
with $(\bX\bbeta_1(P),\bV_0(P))=(\bmu,\bSigma)$.
Let $\bbeta_1(P_{\epsilon,\bs_0})$ and $\bbeta_1(P)$ be a solution to the minimization problem
in Definition~\ref{def:MM-functional general} at $P_{\epsilon,\bs_0}$ and $P$, respectively,
and suppose that
$(\bbeta_1(P_{\epsilon,\bs_0}),\bV_0(P_{\epsilon,\bs_0}))\to (\bbeta_1(P),\bV_0(P))$, as $\epsilon\downarrow0$.
Suppose that $\rho_1:\R\to[0,\infty)$ either satisfies (R-BND), (R-CD1) and (R-CD2),
or satisfies (R-UNB) and (R-CD2), such that~$u_1(s)$ and $u_1'(s)s$ are bounded.
Let
\begin{equation}
\label{def:alpha}
\alpha_1
=
\mathbb{E}_{\mathbf{0},\bI_k}
\left[
\left(1-\frac{1}{k}\right)
\frac{\rho_1'(\|\bz\|)}{\|\bz\|}
+
\frac1k
\rho_1''(\|\bz\|)
\right],
\end{equation}
and suppose that $\mathbb{E}_{0,\bI_k}
\left[
\rho_1''(\|\by\|)
\right]>0$.
If $\bX$ has full rank with probability one, then
for
$\bs_0=(\by_0,\bX_0)$ we have
\[
\mathrm{IF}(\bs_0,\bbeta_1,P)
=
\frac{u_1(d_0)}{\alpha_1}
\Big(
\E\left[\bX^T\bSigma^{-1}\bX\right]
\Big)^{-1}
\bX_0^T\bSigma^{-1}(\by_0-\bX_0\bbeta_1(P)),
\]
where $d_0=d(\bs_0,\bbeta_1(P),\bSigma)$, as defined in~\eqref{def:mahalanobis}.
\end{corollary}
\begin{proof}
Consider $\partial \Lambda_{1}/\partial\bbeta$.
We have
\begin{equation}
\label{eq:derivative psi1 wrt beta}
\frac{\partial\Psi_1(\bs,\bbeta,\bV)}{\partial\bbeta}
=
-
\frac{u_1'(d)}{d}
\bX^T\bV^{-1}(\by-\bX\bbeta)
(\by-\bX\bbeta)^T
\bV^{-1}\bX
-
u_1(d)\bX^T\bV^{-1}\bX,
\end{equation}
where $d=d(\bs,\bbeta,\bV)$, as defined in~\eqref{def:mahalanobis}.
Since $u_1(s)$ and $u_1'(s)s=\rho_1''(s)-u_1(s)$ are bounded,
similar to the proof of Lemma~B.3 in~\cite{lopuhaa-gares-ruizgazenARXIVE2022}, it follows that
for $(\bbeta,\bV)$ in the neighborhood~$N$
of~$(\bbeta_1(P),\bV_0(P))$, it holds that
\begin{equation}
\label{eq:derivative Lambda1}
\frac{\partial\Lambda_1(\bbeta,\bV)}{\partial\bbeta}
=
\int
\frac{\partial\Psi_1(\bs,\bbeta,\bV)}{\partial \bbeta}\,\dd P(\bs),
\end{equation}
and that $\partial\Lambda_1/\partial\bbeta$ is continuous at $(\bbeta_1(P),\bV_0(P))$.
Then, similar to the first part of the proof of Lemma~2 in~\cite{lopuhaa-gares-ruizgazenARXIVE2022},
it can be shown that
\[
\bD_1
=
\frac{\partial\Lambda_1(\bbeta_1(P),\bV_0(P))}{\partial \bbeta}
=
-\alpha_1
\mathbb{E}\left[\mathbf{X}^T\bSigma^{-1}\mathbf{X}\right].
\]
The corollary then follows from Theorem~\ref{th:IF MM}.
\end{proof}
Since $u_1(s)s$ is bounded, the influence function is uniformly bounded in $\by_0$, but not in~$\bX_0$.
This illustrates the phenomenon in linear regression that leverage
points can have a large effect on the regression estimator.

The expression found in Corollary~\ref{cor:IF MM} is the same as the one found for the
regression S-functional in~\cite{lopuhaa-gares-ruizgazenARXIVE2022}
defined with the function $\rho_1$ (see their Corollary~5).
For the multivariate location-scale model, for which $\bX=\bI_k$,
Theorem~\ref{th:IF MM} coincides with Theorem~4.2 in~\cite{lopuhaa1992highly}
and Corollary~\ref{cor:IF MM} matches with the results found in~\cite{SalibianBarrera-VanAelst-Willems2006}.
Furthermore, the expressions found in Theorem~\ref{th:IF MM} and Corollary~\ref{cor:IF MM}
are similar to the ones obtained for the regression MM-functionals in~\cite{yohai1987} and~\cite{kudraszow-maronna2011}, respectively.
For the influence function of MM-functionals in linear mixed effects models, nothing seems to be available yet.
For model~\eqref{def:linear mixed effects model Copt},
the expression for the influence function now follows from Theorem~\ref{th:IF MM}.
For the special case of this model with multivariate normal errors, as considered in~\cite{copt&heritier2007},
the expression for the influence function can be obtained from Corollary~\ref{cor:IF MM}.

\begin{remark}
\label{rem:IF multivariate linear regression}
The multivariate linear regression model~\eqref{def:multivariate linear regression model} is obtained from~\eqref{def:model}
by taking~$\bX=\bx^T\otimes\bI_k$ and $\bbeta=\vc(\bB^T)$.
For this model, the expression in Corollary~\ref{cor:IF MM} for~$\bbeta_1(P)=\mathbf{0}$ is similar,
but slightly different from the one in Theorem~4 in~\cite{kudraszow-maronna2011}.
It seems that Theorem~4 in~\cite{kudraszow-maronna2011} contains some typos
(as confirmed by personal communication).
When~$\bT_1$ is the regression MM-functional considered in~\cite{kudraszow-maronna2011}, then in our notation
$\bbeta_1=\vc(\bT_1^T)$.
When~$P$ is such that $\by\mid\bx$ has an elliptically contoured density
with parameters~$\bmu_0=\bB_0^T\bx$ and $\bSigma_0$,
the correct expression for the influence function of $\bT_1^T$ should be
\[
\text{\rm IF}(\by_0,\bx_0;\bT_1^T,P)
=
\frac{1}{\alpha_1}
u_1\left(
\sqrt{(\by_0-\bB_0^T\bx_0)\bSigma_0^{-1}(\by_0-\bB_0^T\bx_0)}
\right)
(\by_0-\bB_0^T\bx_0)\bx_0^T
\left(
\mathbb{E}\left[\bx\bx^T\right]
\right)^{-1},
\]
with $\alpha_1$ defined in~\eqref{def:alpha} and $u_1(s)=\rho_1'(s)/s$.
\end{remark}

\section{Asymptotic Normality}
\label{sec:asymptotic normality}
Corollary~\ref{cor:consistency MM-estimator bounded} and Theorem~\ref{th:consistency MM unbounded rho1}
provide conditions under which $\bbeta_{1,n}\to\bbeta_1(P)$, with probability one, for $\rho_1$ satisfying
either (R-BND) or (R-UNB).
The next theorem establishes asymptotic normality for $\bbeta_{1,n}$ defined
with either a bounded or an unbounded function $\rho_1$.

\begin{theorem}
\label{th:asymp normal MM}
Suppose that $\rho_1:\R\to[0,\infty)$ either satisfies (R-BND) and (R-CD1), or satisfies (R-UNB).
Suppose that $u_1$ is of bounded variation and let $\E\|\mathbf{s}\|^2<\infty$.
Let $\bbeta_{1,n}$ and~$\bbeta_1(P)$ be solutions to the minimization problems
in Definitions~\ref{def:MM-estimator general} and~\ref{def:MM-functional general}, respectively.
Suppose that $(\bbeta_{1,n},\bV_{0,n})\to (\bbeta_1(P),\bV_0(P))$, in probability,
and suppose that~$\bbeta_1(P)$ is a point of symmetry of $P$.
Suppose that $\Lambda_1$, as defined in~\eqref{def:Lambda beta MM general},
has a partial derivative~$\partial\Lambda_1/\partial\bbeta$ that is continuous
at $(\bbeta_1(P),\bV_0(P))$ and that
$\mathbf{D}_1=(\partial\Lambda_1/\partial\bbeta)(\bbeta_1(P),\bV_0(P))$ is non-singular.
Then $\sqrt{n}(\bbeta_{1,n}-\bbeta_1(P))$ is asymptotically normal with mean zero and
covariance matrix
$\bD_1^{-1}
\E
\left[
\Psi_1(\mathbf{s},\bbeta_1(P),\bV_0(P))
\Psi_1(\mathbf{s},\bbeta_1(P),\bV_0(P))^T
\right]
\bD_1^{-1}$.
\end{theorem}
\begin{proof}
Recall that the estimator can be written as $\bbeta_{1,n}=\bbeta_1(\mathbb{P}_n)$.
This means that it satisfies~\eqref{eq:score psi MM Yohai} for $P=\mathbb{P}_n$:
\begin{equation}
\label{eq:score psi MM estimator Yohai}
\int
\Psi_1(\bs,\bbeta_{1,n},\bV_{0,n})
\,\dd
\mathbb{P}_n(\bs)
=
\mathbf{0},
\end{equation}
where $\Psi_1$ is defined in~\eqref{def:psi beta MM general}.
Writing $\bxi_n=(\bbeta_{1,n},\bV_{0,n})$ and $\bxi_P=\bxi(P)$, we decompose~\eqref{eq:score psi MM estimator Yohai} as follows
\begin{equation}
\label{eq:decomposition MM estimator Yohai}
\begin{split}
\mathbf{0}=
\int \Psi_1(\mathbf{s},\bxi_n)\,\dd P(\mathbf{s})
&+
\int \Psi_1(\mathbf{s},\bxi_P)\,\dd (\mathbb{P}_n-P)(\mathbf{s})\\
&+
\int
\left(
\Psi_1(\mathbf{s},\bxi_n)-\Psi_1(\mathbf{s},\bxi_P)
\right)
\,\dd (\mathbb{P}_n-P)(\mathbf{s}).
\end{split}
\end{equation}
According to~Lemma~B.8 in~\cite{lopuhaa-gares-ruizgazenARXIVE2022}, the third term is of the order $o_P(1/\sqrt{n})$,
whereas according to the central limit theorem the second term is of the order $O_P(1/\sqrt{n})$.
This means we can write
$\mathbf{0}=\Lambda_1(\bxi_n)+O_P(1/\sqrt{n})$,
where
\[
\begin{split}
\Lambda_1(\bxi_n)
&=
\Lambda_1(\bbeta_1(P),\bV_{0,n})
+
\frac{\partial \Lambda_1}{\partial \bbeta}
(\bbeta_1(P),\bV_{0,n})
(\bbeta_{1,n}-\bbeta_1(P))
+
o(\|\bbeta_{1,n}-\bbeta_1(P)\|)\\
&=
\Lambda_1(\bbeta_1(P),\bV_{0,n})
+
\Big(
\bD_1+o_P(1)
\Big)
\Big(\bbeta_{1,n}-\bbeta_1(P)\Big)
+
o(\|\bbeta_{1,n}-\bbeta_1(P)\|).
\end{split}
\]
Since $\bbeta_1(P)$ is a point of symmetry of $P$, according to~\eqref{eq:point of symmetry} it holds that
$\Lambda_1(\bbeta_1(P),\bV_{0,n})=\mathbf{0}$.
It follows that
\[
\mathbf{0}=
\bD_1
(\bbeta_{1,n}-\bbeta_1(P))
+
o(\|\bbeta_{1,n}-\bbeta_1(P)\|)
+
O_P(1/\sqrt{n}).
\]
Because $\bD_1$ is non-singular,
this implies
$\bbeta_{1,n}-\bbeta_1(P)=O(1/\sqrt{n})$.
After inserting this in~\eqref{eq:decomposition MM estimator Yohai}, it follows that
\[
\mathbf{0}
=
\bD_1
\Big(\bbeta_{1,n}-\bbeta_1(P)\Big)
+
\int \Psi_1(\mathbf{s},\bxi_P)\,\dd (\mathbb{P}_n-P)(\mathbf{s})\\
+o_P(1/\sqrt{n}).
\]
We conclude
\[
\sqrt{n}(\bbeta_{1,n}-\bbeta_1(P))
=
-\bD_1^{-1}
\sqrt{n}\int \Psi_1(\mathbf{s},\bxi_P)\,\dd (\mathbb{P}_n-P)(\mathbf{s})\\
+
o_P(1).
\]
Since $\E[\Psi_1(\mathbf{s},\bxi_P)]=\Lambda_1(\bxi_P)=\mathbf{0}$, it follows that
\[
\sqrt{n}\int \Psi_1(\mathbf{s},\bxi_P)\,\dd (\mathbb{P}_n-P)(\mathbf{s})
=
\frac{1}{\sqrt{n}}
\sum_{i=1}^{n}
\Psi_1(\mathbf{s}_i,\bxi_P),
\]
converges in distribution to a multivariate normal random vector with mean zero and covariance
$\E
\left[
\Psi_1(\mathbf{s},\bxi_P)\Psi_1(\mathbf{s},\bxi_P)^T
\right]$.
This finishes the proof.
\end{proof}
When $P$ is such that $\by\mid\bX$ has an elliptically contoured density~\eqref{eq:elliptical}
we can obtain a more detailed expression for the asymptotic covariance.
\begin{corollary}
\label{cor:asymp norm MM}
Suppose that $P$ is such that $\by\mid\bX$ has an elliptically contoured
density~$f_{\bmu,\bSigma}$ from~\eqref{eq:elliptical} with parameters $(\bmu,\bSigma)$.
Suppose that $\bbeta_1(P)$ is the unique minimizer of $R_P(\bbeta)$,
such that $\bX\bbeta_1(P)=\bmu$ and suppose that $\bV_0(P)=\bSigma$.
Let $\bbeta_{1,n}$ and $\bbeta_1(P)$ be solutions to the minimization problems
in Definitions~\ref{def:MM-estimator general} and~\ref{def:MM-functional general}, respectively,
and suppose that $(\bbeta_{1,n},\bV_{0,n})\to (\bbeta_1(P),\bV_0(P))$, in probability.
Suppose that $\rho_1:\R\to[0,\infty)$ either satisfies (R-BND), (R-CD1) and (R-CD2)
or satisfies (R-UNB) and (R-CD2), such that $u_1(s)$ is of bounded variation and $u_1'(s)s$ is bounded,
and let $\E\|\mathbf{s}\|^2<\infty$.
Let~$\alpha_1$ be defined in~\eqref{def:alpha} and suppose that $\mathbb{E}_{0,\bI_k}
\left[
\rho_1''(\|\by\|)
\right]>0$.
If $\bX$ has full rank with probability one, then
$\sqrt{n}(\bbeta_{1,n}-\bbeta_1(P))$ is asymptotically normal with mean zero and
covariance matrix
\begin{equation}
\label{eq:asymp var beta}
\frac{\E_{\mathbf{0},\bI_k}\left[\rho_1'(\|\bz\|)^2\right]}{k\alpha_1^2}
\left(
\mathbb{E}\left[\mathbf{X}^T\bSigma^{-1}\mathbf{X}\right]
\right)^{-1}.
\end{equation}
\end{corollary}
\begin{proof}
Similar to the proof of Corollary~\ref{cor:IF MM} it can be shown that
\[
\bD_1
=
\frac{\partial\Lambda_1(\bbeta_1(P),\bV_0(P))}{\partial \bbeta}
=
-\alpha_1
\mathbb{E}\left[\mathbf{X}^T\bSigma^{-1}\mathbf{X}\right].
\]
According to Theorem~\ref{th:asymp normal MM},
it follows that $\sqrt{n}(\bbeta_n-\bbeta(P))$ is asymptotically normal with mean zero
and covariance matrix
\[
\frac1{\alpha_1^2}
\left(
\mathbb{E}\left[\mathbf{X}^T\bSigma^{-1}\mathbf{X}\right]
\right)^{-1}
\E
\left[
\Psi_1(\mathbf{s},\bbeta_1(P),\bV_0(P))
\Psi_1(\mathbf{s},\bbeta_1(P),\bV_0(P))^T
\right]
\left(
\mathbb{E}\left[\mathbf{X}^T\bSigma^{-1}\mathbf{X}\right]
\right)^{-1},
\]
where $\Psi_1$ is defined in~\eqref{def:psi beta MM general}.
Similar to the proof of Corollary~6 in~\cite{lopuhaa-gares-ruizgazenARXIVE2022}, we find that
\[
\E
\left[
\Psi_1(\mathbf{s},\bbeta_1(P),\bV_0(P))
\Psi_1(\mathbf{s},\bbeta_1(P),\bV_0(P))^T
\right]
=
\frac{\E_{\mathbf{0},\bI_k}\left[\rho_1'(\|\bz\|)^2\right]}{k\alpha_1^2}
\mathbb{E}\left[\mathbf{X}^T\bSigma^{-1}\mathbf{X}\right].
\]
This proves the corollary.
\end{proof}
The expression found in Corollary~\ref{cor:asymp norm MM} coincides with the one for the
regression S-estimator in~\cite{lopuhaa-gares-ruizgazenARXIVE2022} defined with the function $\rho_1$ (see their Corollary~6).
For the multivariate location-scale model, for which $\bX=\bI_k$,
Theorem~\ref{th:asymp normal MM} for unbounded $\rho_1$, coincides with Theorem~3.2 in~\cite{lopuhaa1992highly}.
The expression for the asymptotic variance in Corollary~\ref{cor:asymp norm MM} for bounded $\rho_1$,
coincides
with the results mentioned at the beginning of Section~2.4 in~\cite{SalibianBarrera-VanAelst-Willems2006}.
Furthermore, the results in Theorem~\ref{th:asymp normal MM} and Corollary~\ref{cor:asymp norm MM} are similar to the ones
obtained in~\cite{yohai1987} and~\cite{kudraszow-maronna2011} for regression MM-estimators in the multiple and multivariate linear regression model,
respectively.

For the linear mixed effects model~\eqref{def:linear mixed effects model Copt} with $\bX_i=\bX$ being the same for each subject
and assuming multivariate normal measurement errors, Theorem~1 in~\cite{copt&heritier2007} provides asymptotic normality
of the regression MM-estimator.
Our Theorem~\ref{th:asymp normal MM} and Corollary~\ref{cor:asymp norm MM} are extensions of this result
to a larger class of linear mixed effects models also allowing error distributions much more general than the multivariate normal.

\begin{remark}
\label{rem:Asymp norm multivariate linear regression}
Similar to Remark~\ref{rem:IF multivariate linear regression}, the expression in Corollary~\ref{cor:asymp norm MM}
for the multivariate linear regression model slightly
differs from the one in Proposition~7 in~\cite{kudraszow-maronna2011}.
It seems that the expression in equation~(6.4) in~\cite{kudraszow-maronna2011} contains a small typo
(as confirmed by personal communication).
When~$\widehat{\bB}_n$ is the regression MM-estimator considered in~\cite{kudraszow-maronna2011},
then in our notation $\bbeta_{1,n}=\vc(\widehat{\bB}_n^T)$.
If $P$ is such that $\by\mid\bx$ has an elliptically contoured density
with parameters~$\bmu_0=\bB_0^T\bx$ and $\bSigma_0$,
the correct expression for the asymptotic variance of $\bbeta_{1,n}=\vc(\widehat{\bB}_n^T)$
should be
\[
\frac{\E_{\mathbf{0},\bI_k}\left[\rho_1'(\|\bz\|)^2\right]}{k\alpha_1^2}
\left(
\left(
\E[\bx\bx^T]
\right)^{-1}
\otimes
\bSigma_0
\right),
\]
with $\alpha_1$ defined in~\eqref{def:alpha}.
\end{remark}
Note that for $\rho_{\text{B}}(s;c)$ and $\rho_{\text{H}}(s;c)$, as defined in~\eqref{def:biweight} and~\eqref{def:huber psi}, respectively,
it holds that they both converge to $s^2/2$, as $c\to\infty$.
This means that the least squares estimators can be obtained as limiting case of the regression M-estimator defined with $\rho_1$
equal to either $\rho_{\text{B}}(s;c_1)$ or $\rho_{\text{H}}(s;c_1)$, for $c_1\to\infty$.
In both cases, the scalar
\[
\lambda=
\frac{\E_{\mathbf{0},\bI_k}\left[\rho_1'(\|\bz\|)^2\right]}{k\alpha_1^2}
\to
\frac{\E_{\mathbf{0},\bI_k}\left[\|\bz\|^2\right]}{k},
\quad
\text{as }c_1\to\infty.
\]
For the multivariate normal $\E_{\mathbf{0},\bI_k}\left[\|\bz\|^2\right]=k$, so that the scalar $1/\lambda$ may serve
as an index for the asymptotic efficiency relative to the least squares estimator in all models that are included in our setup.

When using the (bounded) biweight function from~\eqref{def:biweight}, Table~1 in~\cite{kudraszow-maronna2011} gives
the cut-off values $c_0$ for which the initial estimators $(\bbeta_{0,n},\bV_{0,n})$ in Theorem~\ref{th:BDP MM Yohai}
defined with $\rho_0(s)=\rho_{\text{B}}(s;c_0)$ have breakdown point 0.5.
For the regression M-estimator $\bbeta_{1,n}$ defined with $\rho_1(s)=\rho_{\text{B}}(s;c_1)$,
information on cut-off values $c_1$ and asymptotic relative efficiencies can be found at several places in the literature.
Table 1 in~\cite{lopuhaa1989} provides values of $\lambda$, that correspond to the location S-estimator defined with
$\rho_{\text{B}}(s;c_1)$,
for varying dimensions $k=1,2,10$ and varying breakdown points $\epsilon^*=10\%,20\%,\ldots,50\%$, from which $c_1$ can be determined
from
\[
\frac{6\E_{\Phi}
\left[\rho_{\text{B}}(\|\bz\|;c_1)
\right]}{c_1^2}
=
\epsilon^*,
\]
where the expectation is with respect to the standard multivariate normal distribution.
Table~2 in~\cite{kudraszow-maronna2011} provides values of $c_1$ for the multivariate regression MM-estimator,
for asymptotic efficiencies $1/\lambda=80\%,90\%,0.95\%$
and varying dimensions $k=1,\ldots,5,10$.
Finally, Table~3.1 in~\cite{vanaelst&willems2005} gives asymptotic efficiencies that correspond to the
multivariate regression S-estimator for varying dimensions $k=1,2,3,5,10,30,50$ and breakdown points $\epsilon^*=25\%,50\%$.

When using the (unbounded) Huber function from~\eqref{def:huber psi},
Proposition~\ref{prop:BDP MM unbounded} shows that $\bbeta_{1,n}$
inherits the breakdown point of the initial covariance estimator $\bV_{0,n}$,
as long as all $\bX_i=\bX$ are the same and of full rank.
For example, this applies to the linear mixed effects model considered in~\cite{copt&heritier2007}.
For the regression M-estimator $\bbeta_{1,n}$ defined with $\rho_1(s)=\rho_{\text{H}}(s;c_1)$,
information on cut-off values $c_1$ and asymptotic relative efficiencies for the multivariate location M-estimator can be found
in Table~1 in~\cite{maronna1976}
for varying dimensions $k=2,4,6,10$
and ``winsorizing proportions'' $w=0\%, 10\%,20\%,30\%$, from which $c_1$ can be determined via
\[
P_{\Phi}(\|\bz\|>c_1)=w.
\]
Table~1 in~\cite{lopuhaa1989} provides values of $\lambda$, that correspond to the location S-estimator defined with
$\rho_{\text{H}}(s;c_1)$,
for varying dimensions $k=1,2,10$ and the same values for $w$.

\section{Simulation and data example}
\label{sec:application}
We illustrate the finite sample performance of the MM-estimator by means of a simulation.
To this end we will study the behavior of the estimators for samples generated from a model that is similar to the one in~\cite{copt2006high}:
\begin{equation}
\label{eq:model copt}
\by_i=\bX\bbeta+\gamma_i\bZ+\beps_i,
\quad
i=1,\ldots,n,
\end{equation}
a linear mixed effects model with $\by_i$ in dimension $k=4$ and all subjects with the same design matrix $\bX$ for the fixed effects
$\bbeta=(\beta_1,\beta_2)^T$.
Following the setup in~\cite{copt2006high}, the matrix~$\bX$ is built as follows.
The first column of $\bX$ is taken to be a vector $\mathbf{1}$ of length four consisting of ones.
The four $x$-values in the second column are generated from a standard normal,
and then $\bX$ is rescaled to a new matrix $\bX=[\mathbf{1}\quad \bx]$, such that $\bX^T\bX=4\bI_2$.
For our simulation we used
\[
\bX=\left(
  \begin{array}{rr}
1 & -0.9504967\\
1 & -0.5428346\\
1 & 1.6650521\\
1 & -0.1717207\\
\end{array}
\right).
\]
The random effects $\gamma_i$ are independent $N(0,\sigma_\gamma^2)$ distributed random variables, which are independent from the
measurement error $\beps_i\sim N(\mathbf{0},\sigma_\epsilon^2\bI_4)$.
Finally, $\bZ=\mathbf{1}$, a vector of length four consisting of ones.
This leads to a structured covariance $\bSigma=\sigma_\gamma^2\mathbf{1}\mathbf{1}^T+\sigma_\epsilon^2\bI_4$,
with covariance parameter vector $\btheta=(\theta_1,\theta_2)^T$, where $\theta_1=\sigma_\gamma^2$
and $\theta_2=\sigma_\epsilon^2$.
Following the setup in~\cite{copt2006high}, we set $\beta_1=\beta_2=1$ and $\theta_1=\theta_2=1$.

We investigate the behavior of the MM-estimator and the S-estimator used in the first step.
To this end we generate 500 samples of size $n=200$ according to the model in~\eqref{eq:model copt}
and compute the estimates for each sample.
We also study the behavior under contamination, by replacing some of $\by_i$ by $\by_i^*$
obtained by shifting each coordinate of~$\by_i$ over the same distance.
We consider two distances: 3 and 15 and two contamination percentages: 10\% and 20\% contamination.
We used the S-estimator corresponding to Tukey's bi-weight defined in~\eqref{def:biweight}.
The tuning-constant was chosen to be $c_0=4.097$, which corresponds to asymptotic breakdown point 0.5.
For the MM-estimator we used the bi-weight function with tuning constant $c_1=8.530$,
which corresponds to 99\% efficiency relative to the least squares estimator.

\begin{figure}
  \centering
  \begin{tabular}{ccc}
  \includegraphics[width=0.32\textwidth]{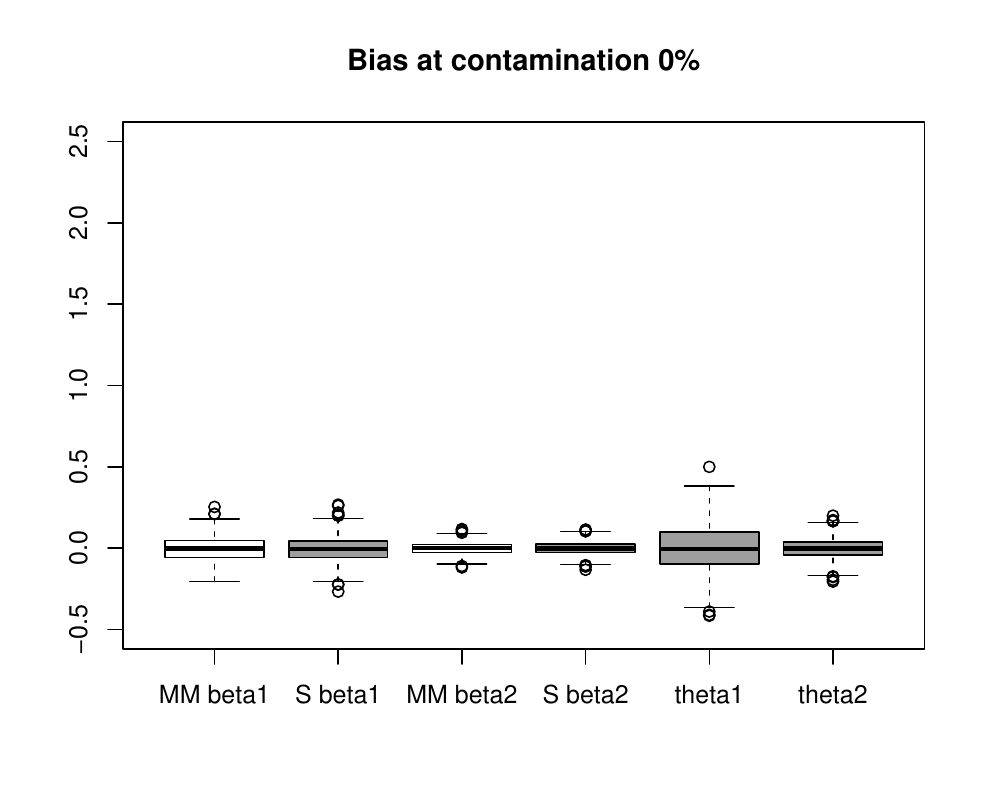}
  &\includegraphics[width=0.32\textwidth]{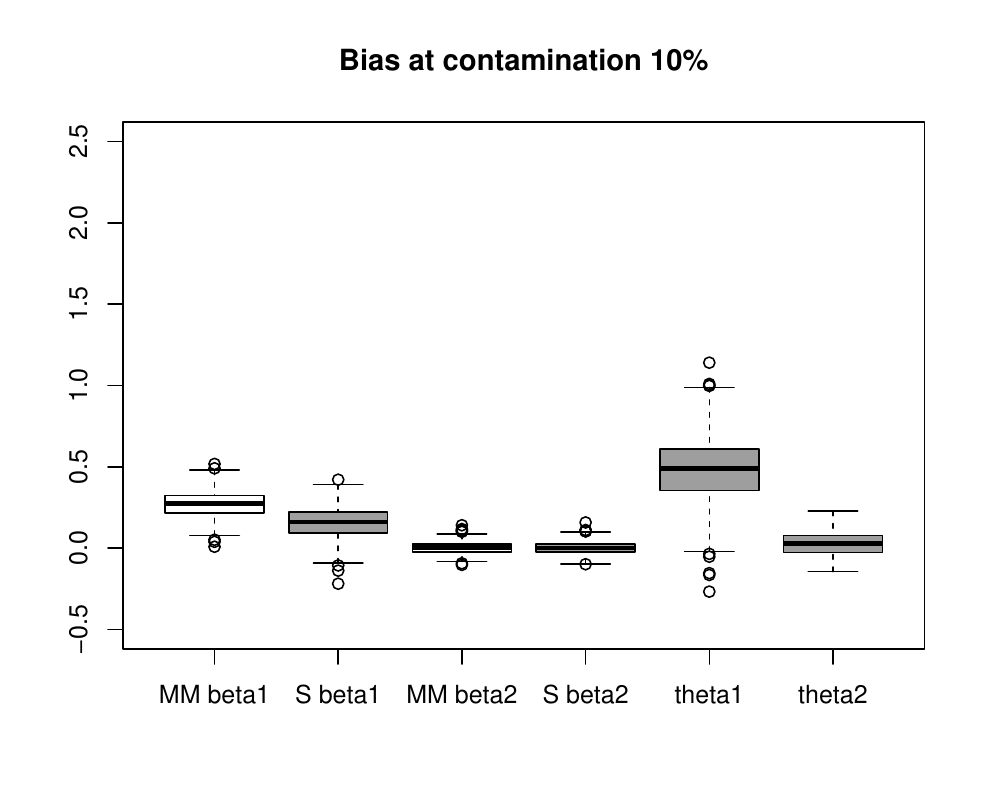}
  &\includegraphics[width=0.32\textwidth]{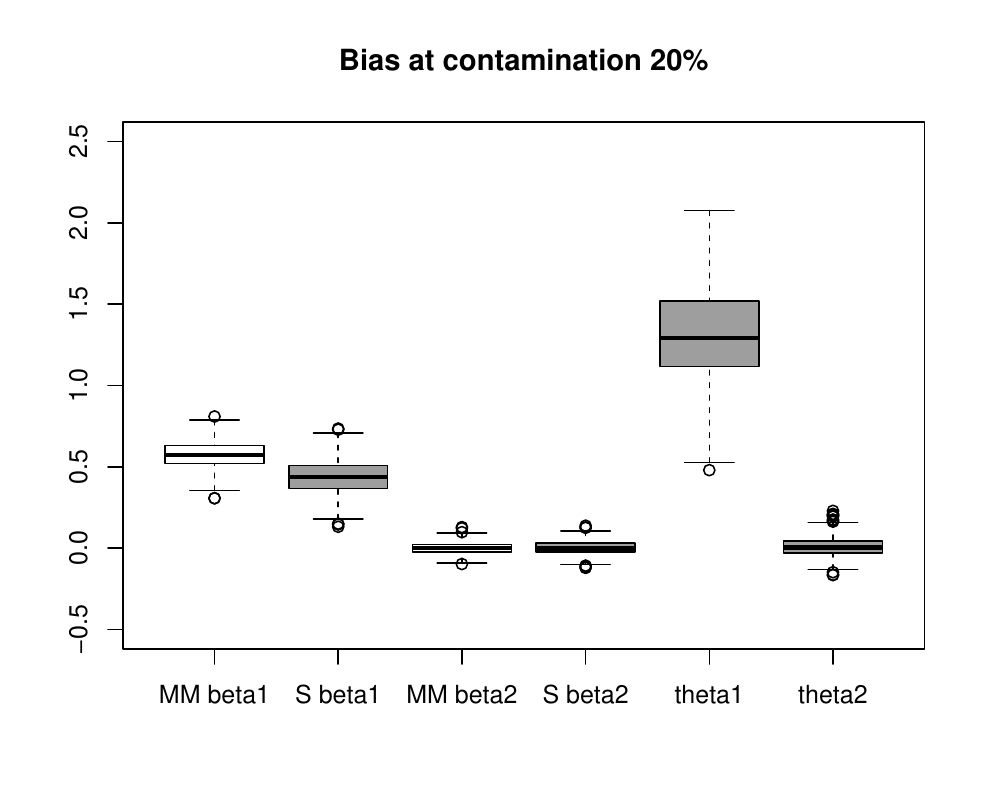}\\
  &\includegraphics[width=0.32\textwidth]{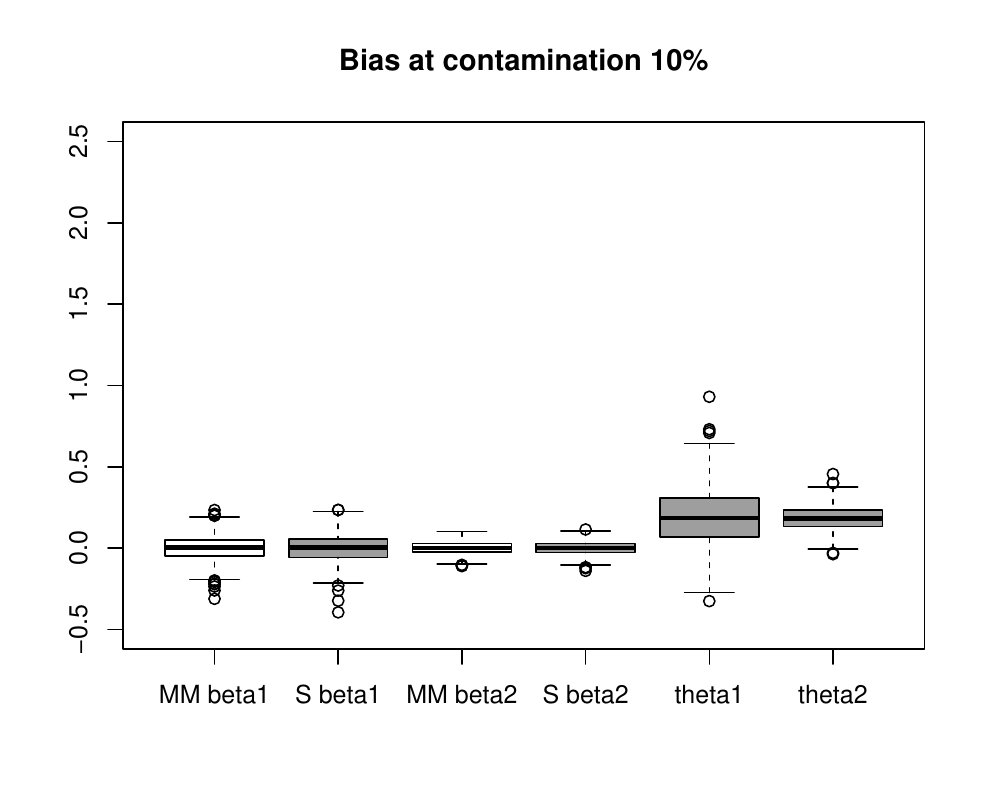}
  &\includegraphics[width=0.32\textwidth]{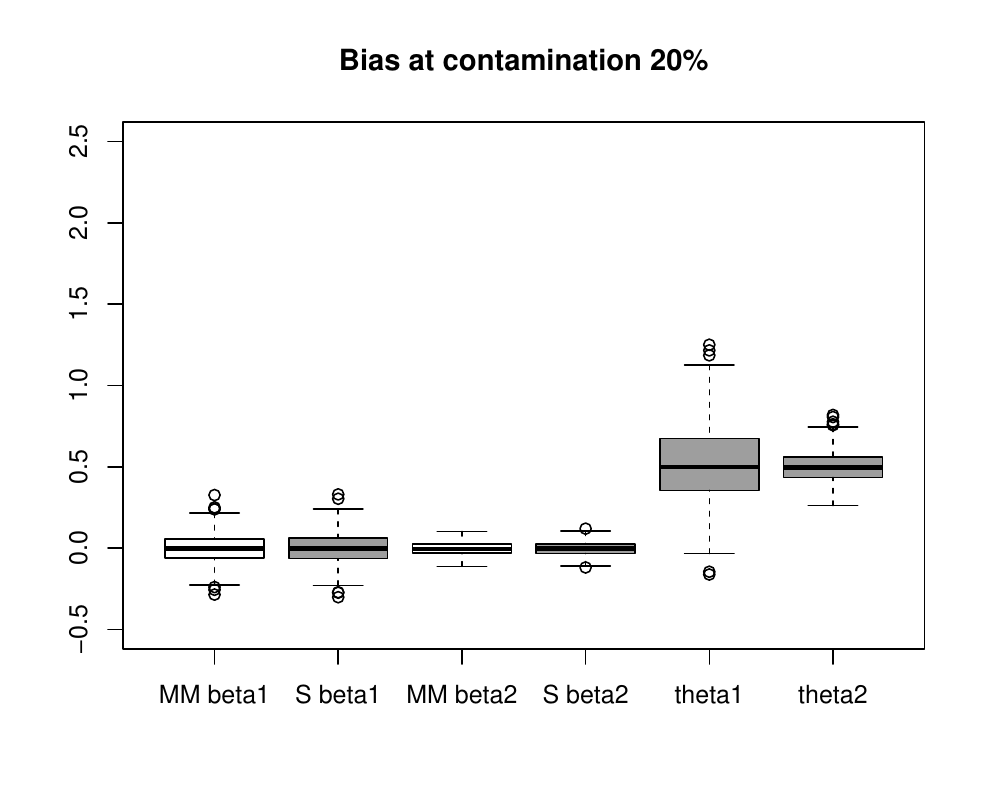}\\
  \end{tabular}
  \caption{MM-estimates and S-estimates at contaminations 0\%, 10\%, and 20\%, with
  shift equal to 3 (top row) and 15 (bottom row).}
  \label{fig:simulation}
\end{figure}

The biases of the estimators are displayed in Figure~\ref{fig:simulation}.
In the case of no contamination, the performance of all estimators is similar;
all the biases are very small.
The fact that the MM-estimator improves the efficiency is illustrated by the fact that
the sample variances of the MM-estimates are 83\% and 82\% of
the sample variances of the S-estimates for $\beta_1$ and $\beta_2$, respectively.
This is comparable to 81\%, which is in both cases the ratio of the two asymptotic variances.
With 10\% contamination by shifting over distance 3, the S-estimator for $\theta_1$ is showing some bias,
as well as the S-estimator for $\beta_1$ to a much lesser extent.
The bias of the S-estimator for $\theta_1$, also seems to affect the MM-estimator for $\beta_1$.
These effects are getting larger when we increase the contamination to 20\%.
When shifting over distance~15, the bias in the estimators for $\beta_1$ disappears,
and is much less for the S-estimator for $\theta_1$, whereas the bias of the S-estimator for $\theta_2$ increases.

This behavior may be explained by the fact that shifting coordinates of contaminated points over distance~3
corresponds to having a difference of about two standard deviations
between the means of the main model and the contaminated one.
In this scenario, the S- and MM-estimators seem to have difficulty separating contaminated points from uncontaminated points.
The S-estimator searches for the smallest cylinder that contains a sufficient number of points
such that it satisfies the constraint in the corresponding minimization problem.
By having contaminated points so close to the bulk of the original data, the S-estimator seems to adapt too much
to the contamination by shifting the intercept of the central axis of the cylinder and stretching
the size of the cylinder along the main principal component, both in the same direction of the shifted points.
This seems to cause the biases in the S-estimators for $\beta_1$ and $\theta_1$.
As a consequence, the MM-estimator is based on points that are standardized by a covariance matrix
that is too large in the direction of the contaminated points.
As such, it will shift even further in the direction of the contaminated points, causing an additional bias in
the MM-estimator for $\beta_1$.
Shifting over distance 15 corresponds to a difference of about ten standard deviations
between the means of the main model and the contaminated one.
In this scenario, the contaminated points are clearly separated from the bulk of the original data.
This seems to make it easier for the S-estimator to find a smallest cylinder $\mathcal{C}(\bbeta,\btheta)$
around the bulk of the original data.
As such, the central axis of the cylinder remains unaffected, resulting in zero bias for the S-estimators for
$\beta_1$ and $\beta_2$.
In order to capture a sufficient number of points to satisfy the constraint of the S-minimization problem,
the S-estimator for the covariance structure is inflated a bit in all directions.
This causes the same amount of bias in the S-estimator for~$\theta_1$ and $\theta_2$, but not as much as in the previous scenario.
Moreover, the MM-estimator is now based on points that are standardized by a covariance matrix that no longer favors a particular direction.
As such, it exhibits zero bias in both 10\% and 20\% contamination.

Finally, we illustrate the performance of the MM-estimator by an application to data from a trial on the treatment of lead-exposed children.
This dataset is discussed in~\cite{fitzmaurice-laird-ware2011} and consists of four repeated measurements of blood lead
levels obtained at baseline (or week~0), week~1, week 4, and week 6 on 100 children who were randomly assigned to chelation treatment
with succimer (a chelation agent) or placebo.
On the basis of a graphical display of the mean response over time, it is suggested in~\cite{fitzmaurice-laird-ware2011} that
a quadratic trend over time seems suitable.
We fitted the following model
\[
y_{ij}=
\beta_0+\beta_1\delta_{i}
+
(\beta_3+\beta_4\delta_i)t_j
+
(\beta_5+\beta_6\delta_i)t_j^2
+
\gamma_{1i}+\gamma_{2i}t_j+\gamma_{3i}t_j^2
+
\epsilon_{ij},
\]
for $i=1,\ldots,100$ and $j=1,\ldots,4$, where
$(t_1,\ldots,t_4)=(0,1,4,6)$ refer to the different weeks,
$y_{ij}$ is the blood lead level (mcg/dL) of subject $i$ obtained at time $t_j$,
and
$\delta_i=0$ if the $i$-th subject is in the placebo group and $\delta_i=1$, otherwise.
The random effects $\bgamma_i=(\gamma_{1i},\gamma_{2i},\gamma_{3i})$, $i=1,\ldots,100$, are assumed to be independent
mean zero normal random vectors with a diagonal covariance matrix consisting of variances $\sigma_{\gamma_1}^2$,
$\sigma_{\gamma_2}^2$ and $\sigma_{\gamma_3}^2$, respectively.
The measurement errors $\beps_i=(\epsilon_{i1},\ldots,\epsilon_{i4})$, $i=1,\ldots,100$, are assumed to be independent
mean zero random vectors with covariance matrix $\sigma_\epsilon^2\bI_4$,
also being independent of the random effects.
In this way we are fitting a balanced linear mixed effects model with unknown parameters $\bbeta=(\beta_1,\ldots,\beta_6)$
and $\btheta=(\sigma_{\gamma_1}^2,\sigma_{\gamma_2}^2,\sigma_{\gamma_3}^2,\sigma_{\epsilon}^2)$.

We estimated $(\bbeta,\btheta)$ by means of least squares and by means of the S-estimator corresponding to
Tukey's bi-weight defined in~\eqref{def:biweight}.
The tuning-constant was chosen to be $c_0=4.097$, which corresponds to asymptotic breakdown point 0.5
and 80\% efficiency relative to the least squares estimator.
For the MM-estimator we used the bi-weight function with tuning constant $c_1=5.810$,
which corresponds to 95\% efficiency relative to the least squares estimator.
The resulting estimates and their standard errors (between brackets) are given in Table~\ref{tab:estimates}.
\begin{table}[t!]
  \centering
\begin{tabular}{crccrccrc}
  \hline
\\[-10pt]
  & \multicolumn{2}{c}{LS} && \multicolumn{2}{c}{S} && \multicolumn{2}{c}{MM}\\
\\[-10pt]
  \cline{2-3} \cline{5-6}\cline{8-9}
\\[-10pt]
$\beta_1$ & 23.973 & (0.899)&& 23.188 & (0.962) && 23.341 & (0.883)\\
$\beta_2$ &  1.996 & (1.272)&&  2.177 & (1.361) &&  2.342 & (1.249)\\
$\beta_3$ & $-7.541$ &(0.591) && $-6.298$ & (0.548) && $-7.047$ & (0.503)\\
$\beta_4$ & 6.624 & (0.835)&& 5.298 & (0.775) && 6.089 & (0.712)\\
$\beta_5$ & 1.196 & (0.097) && 0.985 & (0.085) && 1.105 & (0.078)\\
$\beta_6$ & $-1.104$ & (0.137) && $-0.890$ & (0.120) && $-1.012$ & (0.110 )\\
$\sigma_{\gamma_1}^2$ & 23.253 & && 25.756\\
$\sigma_{\gamma_2}^2$ & $-0.009$ & && $0.575$\\
$\sigma_{\gamma_3}^2$ & 0.004 & && $-0.020$\\
$\sigma_{\epsilon}^2$ & 21.785 & && 14.293\\
\\[-10pt]
\hline
\end{tabular}
  \caption{LS estimates, S-estimates, and MM-estimates for $\beta_1,\ldots,\beta_6$ with their standard errors.}
  \label{tab:estimates}
\end{table}
The standard errors can be computed from~\eqref{eq:asymp var beta},
taking the tuning-constant $c_1$ equal to infinity, $4.097$ and $5.810$, for the LS, S-, and MM-estimator, respectively.
The MM-estimates are more in line with the least squares estimates and have smaller standard errors than those of the S-estimates.
At the same time, the robust standardized residuals computed from the MM-estimates $\bbeta_{1,n}$ and $\bV(\btheta_{0,n})$,
are almost identical to ones computed from the S-estimates $\bbeta_{0,n}$ and $\bV(\btheta_{0,n})$.
This means that the MM-estimator identifies the same observations 40 and 98 as outliers,
whereas the least squares estimator only identifies observation 40 (see also Figure 2 in~\cite{lopuhaa-gares-ruizgazenARXIVE2022}).

\section{Discussion}
\label{sec:discussion}
We have provided a unified approach for constructing estimators of the regression parameter in balanced linear models
with a structured covariance that combine good robustness properties, such as a high breakdown point and an influence function
that is bounded in $\by$,
with high asymptotic efficiency at models with multivariate normal errors.
Our setup is sufficiently flexible to include several specific multivariate statistical models,
including linear mixed effects models, multivariate and multiple linear regression, and multivariate location and scale.
In this way, the theory for these multivariate models can be handled at the same time.

Combining high breakdown and with high efficiency originated with the MM-estimators introduced by~\cite{yohai1987}
for the multiple regression model.
For the models included in our setup, extensions have been developed in two ways.
For example, one approach for multivariate location and scatter
is to start by obtaining an initial high breakdown estimator of the covariance matrix,
and use this to standardize the observations and determine an location M-estimator
from the rescaled observations in such a way that it has high efficiency
and inherits the breakdown point of the covariance estimator, see~\cite{lopuhaa1992highly}.
A similar approach was used in~\cite{copt&heritier2007} for the regression estimator in linear mixed effects models.
Another approach for multivariate location and scatter
is to start by obtaining initial high breakdown estimators of location and of the shape of the scatter matrix,
and use these to determine an auxiliary univariate M-estimator of scale.
The estimators of location and shape are then updated to improve the efficiency,
after which the latter is combined with the M-estimator of scale to determine the final estimator of scatter,
see~\cite{tatsuoka&tyler2000} and~\cite{SalibianBarrera-VanAelst-Willems2006}.
This approach was also used in~\cite{kudraszow-maronna2011} for multivariate linear regression.
Since the main objective for our general setup was to construct an estimator for the regression parameter
that combines high breakdown with high efficiency we extended the somewhat simpler approach used
in~\cite{lopuhaa1992highly,copt&heritier2007}.
The more involved approach used in~\cite{tatsuoka&tyler2000,SalibianBarrera-VanAelst-Willems2006,kudraszow-maronna2011}
does have the advantage that also the efficiency of the covariance shape estimator can be improved.
This approach is beyond the scope of this paper and will be considered in a future manuscript.

Our main interest is to combine high breakdown with high efficiency for linear mixed effects models, for which the theory is far from complete.
The results in the literature are limited to models in which the design matrix for the fixed effects is the same for all subjects,
and only deal with the asymptotic behavior under the assumption of normally distributed errors.
Moreover, no rigorous attention has been paid to whether the estimators actually exist and what their robustness properties are.
Our general setup includes linear mixed effects models,
which allow different subject specific design matrices and
general multivariate distributions that go far beyond the multivariate normal or other elliptically contoured distributions.
We have provided sufficient conditions under which the estimators and corresponding statistical functionals exist.
Furthermore, we have given conditions under which the regression estimator inherits the (high) breakdown point of the covariance estimator in the initial step, and we have derived the expression of the influence function.
Finally, we have established strong consistency and asymptotic normality of the regression estimator under general distributions,
from which more detailed expressions can be determined at the multivariate normal or other elliptically contoured distributions.

Since our setup also includes the multivariate location-scale model, our results also include the ones obtained in~\cite{lopuhaa1992highly}
as a special case.
In addition, our results have extended the ones on existence and breakdown point to location MM-estimators defined by means of a bounded $\rho_1$ in Definition~\ref{def:MM-estimator general}.
Although our simpler approach differs from the one in~\cite{SalibianBarrera-VanAelst-Willems2006},
our results on the breakdown behavior of the location MM-estimator are similar.
Furthermore, our results on the influence function and on the asymptotic behavior of the location MM-estimator are
much more general, but are identical to the ones in~\cite{SalibianBarrera-VanAelst-Willems2006} at the special case of elliptically contoured distributions.
Also the multivariate linear regression model is included in our setup,
although our simpler approach differs from the one in~\cite{kudraszow-maronna2011}.
Again our results on the existence and the breakdown behavior of the regression MM-estimator are similar,
and our results on the influence function and the asymptotic behavior of the regression MM-estimator are more general,
but identical to the ones in~\cite{kudraszow-maronna2011} at the special case of elliptically contoured errors distributions.


\end{document}